\documentclass[11pt]{article}
\usepackage[dvips]{graphics}
\usepackage{epsfig,rotating}
\usepackage{amssymb,amsmath}
\usepackage{latexsym}
\usepackage[usenames]{color}

\newcommand{\be}{\begin{equation}}
\newcommand{\ee}{\end{equation}}
\newcommand{\beaa}{\begin{eqnarray*}}
\newcommand{\eeaa}{\end{eqnarray*}}
\newcommand{\bea}{\begin{eqnarray}}
\newcommand{\eea}{\end{eqnarray}}
\newcommand{\lbl}{\label}

\newcommand{\mn}{|\!|\!|}

\newcommand{\hf}{{1 \over 2}}

\newtheorem{theorem}{ \noindent T{\footnotesize HEOREM}}
\newtheorem{prop}{ \noindent P{\footnotesize ROPOSITION}}[section]
\newtheorem{lemma}{ \noindent L{\footnotesize EMMA}}[section]
\newtheorem{coro}{ \noindent C{\footnotesize OROLLARY}}[section]
\newtheorem{remark}{ \noindent R{\footnotesize EMARK}}[section]

\newcommand{\RR}{\mathbb{R}}
\newcommand{\goto}{\rightarrow}
\begin{document}

\title{Limiting Laws of Coherence of Random Matrices with
  Applications to Testing Covariance Structure
  and Construction of Compressed Sensing Matrices }
\author{Tony Cai$^{1}$ and Tiefeng Jiang$^{2}$\\
University of Pennsylvania and University of Minnesota}

\date{}
\maketitle

\footnotetext[1]{ Statistics Department, The Wharton School, University of
  Pennsylvania,  Philadelphia, PA 19104, \newline  \indent \ \
tcai@wharton.upenn.edu. The research of Tony Cai was supported in part
by NSF FRG Grant  \newline  \indent \ \
DMS-0854973.}
\footnotetext[2]{School of Statistics, University of Minnesota, 224 Church
Street, MN55455, tjiang@stat.umn.edu. \newline  \indent \ \
The research of Tiefeng Jiang was
supported in part by NSF FRG Grant DMS-0449365.}

\begin{abstract}
Testing covariance structure is of significant interest in many areas of
statistical analysis and construction of compressed sensing
matrices is an important problem in signal processing. Motivated by
these applications, we study in this paper the
limiting laws of the coherence of an $n\times p$ random matrix in the
high-dimensional setting where $p$ can be much larger than $n$.
Both the law of large numbers and the limiting distribution are derived.
We then consider testing the bandedness of the covariance matrix of
a high dimensional Gaussian distribution which includes testing for
independence as a special case. The limiting laws of the coherence of
the data matrix play a critical role in the construction of the test.
We also apply the asymptotic results to the construction of compressed
sensing matrices.
\end{abstract}

\noindent \textbf{Keywords:\/}  Chen-Stein method, coherence,
compressed sensing matrix, covariance structure, law of large numbers,
limiting distribution, maxima, moderate deviations, mutual incoherence
property, random matrix, sample correlation matrix.

\noindent\textbf{AMS 2000 Subject Classification: \/} Primary 62H12, 60F05;
secondary 60F15, 62H10.

\newpage

\textbf{\large Contents}

{\small
\ref{intro}. Introduction

\ \ \ \ \ref{LLCR}. Limiting Laws of the Coherence of a Random Matrix

\ \ \ \ \ref{TCS}.  Testing Covariance Structure

\ \ \ \ \ref{CCSM}. Construction of Compressed Sensing Matrices


\ \ \ \ \ref{OP}. Organization of the Paper

\ref{limit.sec}. Limiting Laws of Coherence of Random Matrix

\ \ \ \ \ref{iid.sec}. The i.i.d. Case

\ \ \ \ \ref{dependent.sec}. The Dependent Case

\ref{testing.sec}. Testing the Covariance Structure

\ref{CS.sec}. Construction of Compressed Sensing Matrices

\ref{discussion.sec}. Discussion and Comparison with Related Results

\ref{proof.sec} Proofs

\ \ \ \ \ref{Technical}. Technical Tools

\ \ \ \ \ref{LLN}. Proofs of Theorems \ref{army} and \ref{morethan}

\ \ \ \ \ref{Limiting}. Proof of Theorem \ref{birthday}

\ \ \ \ \ref{substantiald}. Proof of Theorem \ref{youngth}

\ref{ladies} Appendix
}

\section{Introduction}
\lbl{intro}
\setcounter{equation}{0}

Random matrix theory has been proved to be a powerful tool in a wide
range of fields including statistics, high-energy physics, electrical
engineering and number theory. Traditionally the primary focus is on the
spectral analysis  of eigenvalues and eigenvectors. See, for example,
Johnstone (2001 and 2008), Bai, Miao and Pan (2007), and Jiang
(2004b).  For general background on the random matrix theory, see, for example, Bai and
Silverstein (2009) and Anderson, Guionnet, and Zeitouni (2009).

In statistics, the random matrix theory is particularly useful for inference of
high-dimensional data which is becoming increasingly available in many areas
of scientific investigations. In these applications, the dimension $p$
can be much larger than the sample size $n$. In such a setting
classical statistical methods and results based on fixed $p$ and large
$n$ are no longer applicable. Examples include high-dimensional
regression, hypothesis testing concerning high-dimensional parameters,
and inference on large covariance matrices.  See, for example, Candes and
Tao (2007), Cai, Wang and Xu (2010a), Bai and Saranadasa (1996), Bai,
Jiang, Yao and Zheng (2009), and Cai, Zhang and Zhou (2010).

In the present paper  we study the limiting laws of the coherence of an
$n\times p$ random matrix, which is defined to be the largest
magnitude of the off-diagonal entries of the sample correlation
matrix generated from the $n\times p$ random matrix. We are especially
interested in the case where $p\gg n$.
This is a problem of independent interest. Moreover, we are
particularly interested in the applications of the results to
testing the covariance structure of a high-dimensional Gaussian variable and the construction of compressed sensing matrices.
These three problems are important in their respective fields, one in
random matrix theory, one in statistics and one in signal processing.
The latter two problems are seemingly unrelated at first sight, but as we
shall see later they can both be attacked through the use of the
limiting laws of the coherence of random matrices.

\subsection{Limiting Laws of the Coherence of a Random Matrix}\lbl{LLCR}

Let $X_n=(x_{ij})$ be an $n\times p$ random matrix where the entries
$x_{ij}$ are i.i.d. real random variables with mean
$\mu$ and variance $\sigma^2>0$. Let $x_1, x_2, \cdots, x_p$ be the $p$
columns of $X_n$. The sample correlation matrix $\Gamma_n$ is defined by
$\Gamma_n:=(\rho_{ij})$ with
\bea\lbl{corr}
\rho_{ij}=\frac{(x_i-\bar{x}_i)^T(x_j-\bar{x}_j)}{\|x_i-\bar{x}_i\|\cdot \|x_j-\bar{x}_j\|},\ \ 1\leq i, j\leq p
\eea
where $\bar{x}_k=(1/n)\sum_{i=1}^nx_{ik}$ and $\|\cdot \|$ is the
usual Euclidean norm in $\mathbb{R}^n$. Here we write $x_i-\bar{x}_i$
for $x_i-\bar{x}_ie,$ where $e=(1,1,\cdots, 1)^T \in \mathbb{R}^n.$
In certain applications such as construction of compressed sensing
matrices, the mean $\mu$ of the random entries $x_{ij}$ is known (typically
$\mu=0$) and the sample
correlation matrix is then defined to be $\tilde{\Gamma}_n:=(\tilde{\rho}_{ij})$ with
\be\lbl{corr'}
\tilde{\rho}_{ij}=\frac{(x_i-\mu)^T(x_j-\mu)}{\|x_i-\mu\|\cdot \|x_j-\mu\|},\
\ 1\leq i, j\leq p.
\ee

One of the main objects of interest in the present paper is the largest
magnitude of the off-diagonal entries of the sample correlation
matrix,
\bea\lbl{mount}
L_n=\max_{1\leq i <  j \leq p}|\rho_{ij}| \ \ \mbox{and}\ \
\tilde{L}_n=\max_{1\leq i <  j \leq p}|\tilde{\rho}_{ij}|.
\eea
In the compressed sensing literature, the quantity $\tilde{L}_n$ is called the
{\it coherence} of the matrix $X_n$. A matrix is incoherent when
$\tilde{L}_n$ is small. See, for example, Donoho, Elad
and Temlyakov (2006). With slight abuse of terminology,
in this paper we shall call both $L_n$ and $\tilde{L}_n$ {\it coherence} of
the random matrix $X_n$, the former for the case $\mu$ is unknown and
the latter for the case $\mu$ is known. The first goal of the present
paper  is to derive
the limiting laws of the coherence in the high dimensional setting.

In the case where $p$ and $n$ are comparable, i.e., $n/p \goto \gamma
\in (0,\infty)$, asymptotic properties of the coherence $L_n$ of
random matrix $X_n$ have been considered by Jiang (2004a), Zhou (2007), Liu, Lin and Shao (2008),
and Li, Liu and Rosalsky (2009). In this paper we focus on the high
dimensional case where $p$ can be as large as $e^{n^\beta}$ for some
$0<\beta<1$. This is a case of special interest for the applications considered later.

The results given in Section \ref{limit.sec} show that under
regularity conditions,
\[
\sqrt{n/\log p}\, L_n \stackrel{P}{\to} 2 \quad \mbox{as} \;\; n\to\infty
\]
where $\stackrel{P}{\to}$ denotes convergence in
probability. Here and throughout the paper the $\log$ is the natural logarithm $\log_e.$ Furthermore, it is shown that
$nL_n^2-4\log p + \log\log p$ converges weakly to an extreme
distribution of type I with
distribution function
\[
F(y)=e^{-{1\over \sqrt{8\pi}}e^{-y/2}}, \;\; y\in \mathbb{R}.
\]
Same results hold for $\tilde{L}_n$. In contrast to the known results in the
literature, here the dimension $p$ can be much larger
than $n$. In the special cases where $x_{ij}$ are either bounded or
normally distributed, the results hold as long as  $\log p = o(n^{1/3})$.

In addition, motivated by application to testing covariance structure,
we also consider the case where the entries of random
matrix $X_n$ are correlated. More specifically, let
$X_n=(x_{ij})_{1\leq i \leq n, 1\leq j \leq p},$ where the $n$ rows
are i.i.d. random vectors with distribution $N_p(\mu, \Sigma)$.
For a given integer $\tau\geq 1$ (which can depend on $n$ or $p$), it is of
interest in applications to test the hypothesis that the covariance
matrix $\Sigma$ is banded, that is,
\begin{eqnarray}\lbl{Hypo100}
H_0: \sigma_{ij}=0\ \mbox{for all } |i-j|\geq \tau.
\end{eqnarray}
Analogous to the definition of $L_n$ and $\tilde{L}_n$, we define
\be\lbl{lag-k}
L_{n,\tau}=\max_{|i-j|\geq \tau}|\rho_{ij}|
\ee
when the mean $\mu$ is assumed to be unknown and define
\be\lbl{lag-k'}
\tilde{L}_{n,\tau}=\max_{|i-j|\geq \tau}|\tilde{\rho}_{ij}|
\ee
when the mean $\mu=(\mu_1, \mu_2, ..., \mu_p)$ is assumed to be
known. In the latter case $\tilde{\rho}_{i,j}$ is defined to be
\be\lbl{corr''}
\tilde{\rho}_{ij}=\frac{(x_i-\mu_i)^T(x_j-\mu_j)}{\|x_i-\mu_i\|\cdot \|x_j-\mu_j\|},\
\ 1\leq i, j\leq p.
\ee
We shall derive in Section \ref{limit.sec} the limiting distribution
of $L_{n,\tau}$ and $\tilde{L}_{n,\tau}$ under the null hypothesis $H_0$ and discuss its
application in Section \ref{testing.sec}.
The study for this case is considerably more difficult technically than that for the i.i.d. case.

\subsection{Testing Covariance Structure}\lbl{TCS}

Covariance matrices play a critical role in many areas of
statistical inference. Important examples include principal
component analysis, regression analysis, linear and quadratic
discriminant analysis, and graphical models. In the classical setting
of low dimension and large sample size, many methods have been
developed for estimating covariance matrices as well as testing
specific patterns of covariance matrices.
In particular testing for independence in the Gaussian case is of
special interest because many statistical procedures are built upon
the assumptions of independence and normality of the observations.

To be more specific, suppose we observe independent and identically
distributed $p$-variate random variables
$\mathbf{Y}_{1},\ldots ,\mathbf{Y}_{n}$ with mean $\mu=\mu_{p\times 1}$,
covariance matrix $\Sigma=\Sigma_{p\times p}$ and correlation matrix
$R=R_{p\times p}$.
In the setting where the dimension $p$ and the sample size $n$ are
comparable, i.e., $n/p \goto \gamma \in (0,\infty)$, testing of the
hypotheses $H_0: \Sigma = I$ versus $H_a: \Sigma \neq I$, assuming
$\mu=0$, has been considered by Johnstone (2001) in the Gaussian case
and  by P\'{e}ch\'{e} (2009) in the more general case where the
distribution is assumed to be sub-Gaussian and where the ratio $p/n$
can converge to either a positive number $\gamma$,
0 or $\infty$.
The test statistic is based on the largest eigenvalue of the sample
covariance matrix and relies on the important results in their papers
that the largest eigenvalue of the sample covariance matrix
follows the Tracy-Widom distribution asymptotically.

The hypothesis $H_0: \Sigma = I$ is too restrictive for many
applications. An arguably more practically important problem is
testing for independence in the Gaussian
case. That is, one wishes to test the
hypothesis $H_0: \Sigma$ is diagonal against the hypothesis $H_a:
\Sigma$ is not diagonal, or equivalently in terms of the correlation
matrix $R$, one wishes to test $H_0: R = I$ versus
$H_a: R \neq I$. Tests based on the largest eigenvalue of the
sample covariance matrix cannot be easily modified for testing these
hypotheses.

In this paper, we consider testing more general hypotheses on the
covariance structure of a high dimensional Gaussian distribution which
includes testing for independence as a special case. More specifically,
we consider testing the hypothesis that $\Sigma$
is banded with a given bandwidth $\tau$ (which may depend on $n$ or $p$), i.e., the variables have nonzero
correlations only up to lag $\tau$. In other words, for a given integer
$\tau\ge 1$, we wish to test the hypothesis $H_0$: $\sigma_{i,j} = 0$ for
all $|i-j|\ge \tau$.
This problem arises, for example, in econometrics when testing certain
economic theories and in time series analysis. See Andrews (1991),
Ligeralde and Brown (1995) and references therein. The special case of
$\tau=1$ corresponds to testing for independence.
We shall show that the limiting laws of $L_{n,\tau}$ developed in
the present paper can be readily  applied to construct a convenient test for
the bandedness of the covariance matrix. In the special case of $\tau=1$,
the limiting laws of the coherence of the data matrix $\mathbf{Y}$ play a
critical role in the construction of the test.

\subsection{Construction of Compressed Sensing Matrices}\lbl{CCSM}

In addition to testing the covariance structure, another important
application of our results on the limiting laws of the
coherence of a random matrix is to the construction of compressed
sensing  matrices.
Compressed sensing is a fast developing field which provides a novel
and efficient data acquisition technique that enables accurate
reconstruction of highly undersampled sparse signals. See, for
example, Donoho (2006a).
It has a wide range of applications including signal processing,
medical imaging, and seismology. In addition, the development
of the compressed sensing theory also provides crucial insights into
high dimensional regression in statistics. See, e.g., Candes and Tao
(2007), Bickel, Ritov and Tsybakov (2009),  and Candes and Plan (2009).

One of the main goals of compressed sensing is to construct measurement
matrices $X_{n\times p}$, with the number of measurements $n$ as small
as possible relative to $p$, such that for any $k$-sparse signal
$\beta \in \RR^p$,
one can recover $\beta$ exactly from linear measurements $y =
X\beta$  using a computationally efficient recovery algorithm. In
compressed sensing it is typical that $p\gg n$, for example, $p$ can
be order $e^{n^\beta}$ for some $0<\beta<1$. In fact, the goal is
often to make $p$ as large as possible relative to $n$. It is now well
understood that the method of $\ell_1$ minimization provides an
effective way for reconstructing a sparse signal in many settings. In
order for a recovery algorithm such as $\ell_1$ minimization to work
well, the measurement
matrices $X_{n\times p}$ must satisfy certain conditions. Two commonly
used conditions are the so called restricted isometry property (RIP)
and  {\it mutual incoherence property} (MIP). Roughly speaking, the RIP
requires subsets of certain cardinality of the columns of $X$ to be
close to an orthonormal system and the MIP requires the pairwise
correlations among the column vectors of $X$ to be small.
See Candes and Tao (2005), Donoho, Elad and Temlyakov (2006) and Cai,
Wang and Xu (2010a, b).
It is well known that construction of large deterministic measurement
matrices that satisfy either the RIP or MIP  is
difficult. Instead, random matrices are commonly used. Matrices
generated by certain random processes have been shown to satisfy the
RIP conditions with high probability. See, e.g.,  Baraniuk, et. al. (2008).
A major technical tool used there is the Johnson-Lindenstrauss lemma.
Here we focus on the MIP.

The MIP condition can be easily explained. It was  first shown by
Donoho and Huo (2001), in the setting where $X$ is a concatenation of
two square orthogonal matrices, that the condition
\be
\label{sharp.condition}
(2k-1)\tilde{L}_n < 1
\ee
ensures the exact recovery of $\beta$ when $\beta$ has at most $k$
nonzero entries (such a signal is called $k$-sparse). This result was
then extended by Fuchs (2004) to general matrices. Cai,
Wang and Xu (2010b) showed that condition (\ref{sharp.condition})
is also sufficient for stable recovery of sparse signal in the noisy
case where $y$ is measured with error. In addition, it was shown that
this condition is sharp in the sense that there exist matrices $X$
such that it is not possible to recover certain $k$-sparse signals
$\beta$ based on $y=X\beta$ when $(2k-1)\tilde{L}_n =1$.

The mutual incoherence property (\ref{sharp.condition}) is  very
desirable. When it is satisfied by the measurement matrix
$X$, the estimator obtained through $\ell_1$ minimization satisfies
near-optimality properties and oracle inequalities.
In addition, the technical analysis is particularly simple. See, for
example,  Cai, Wang and Xu (2010b). Except results on the magnitude
and the limiting distribution of $\tilde L_n$ when the underlying matrix is Haar-invariant and orthogonal by Jiang (2005),   it is, however, unknown in general
how likely a random matrix satisfies the MIP (\ref{sharp.condition})
in the high dimensional setting where $p$ can be as large as $e^{n^\beta}$.
We shall show in Section \ref{CS.sec} that the limiting laws of the
coherence of random matrices given in this paper can readily be
applied to compute the probability that random measurement matrices satisfy
the MIP condition (\ref{sharp.condition}).

\subsection{Organization of the Paper}\lbl{OP}

The rest of the paper is organized as follows. We begin in Section
\ref{limit.sec} by studying the limiting laws of the coherence of a
random matrix in the high-dimensional
setting. Section \ref{testing.sec} considers the problem of testing for
independence and bandedness  in the Gaussian case. The test
statistic is based on the coherence of the data matrix and the
construction of the tests relies heavily on the asymptotic results
developed in Section \ref{limit.sec}. Application to the construction
of compressed sensing matrices is considered in Section \ref{CS.sec}.
Section \ref{discussion.sec} discusses connections and differences of
the our results with other related work.
The main results are proved in Section \ref{proof.sec} and the proofs of
technical lemmas are given in the Appendix.

\section{Limiting Laws of Coherence of Random Matrices}
\label{limit.sec}

In this section, we consider the limiting laws of the coherence of a
random  matrix with i.i.d. entries. In addition, we also consider the
case where each row of the random matrix is drawn independently from a
multivariate Gaussian distribution with banded covariance matrix. In
the latter case we consider the limiting distribution of $L_{n,\tau}$ and
$\tilde{L}_{n,\tau}$ defined in (\ref{lag-k}) and (\ref{lag-k'}).
We then apply the asymptotic results to the testing of the covariance
structure in Section \ref{testing.sec} and the construction of
compressed sensing matrices in Section \ref{CS.sec}.

\subsection{The i.i.d. Case}
\label{iid.sec}

We begin by considering the case for independence where all entries of the
random matrix are independent and identically distributed.
Suppose $\{\xi, \; x_{ij},\ i,j=1,2,\cdots\}$ are
i.i.d. real random variables with mean $\mu$ and
variance $\sigma^2>0$. Let $X_n=(x_{ij})_{1\leq
  i \leq n, 1\leq j \leq p}$ and let $x_1, x_2, \cdots, x_p$ be the
$p$ columns of $X_n.$ Then $X_n=(x_1, x_2, \cdots, x_p).$ Let
$\bar{x}_k=(1/n)\sum_{i=1}^nx_{ik}$ be the sample average of $x_k$.
We write $x_i-\bar{x}_i$ for $x_i-\bar{x}_ie,$ where $e=(1,1,\cdots,
1)^T \in \mathbb{R}^n.$  Define the Pearson correlation coefficient
$\rho_{ij}$  between $x_i$ and $x_j$ as in (\ref{corr}).
Then the {\it sample correlation matrix} generated by $X_n$  is
$\Gamma_n:=(\rho_{ij})$, which is a $p$ by $p$ symmetric matrix with diagonal
entries $\rho_{ii}=1$ for all $1\le i\le p$.
When the mean $\mu$ of the random variables $x_{ij}$ is assumed to be
known, we define the sample correlation matrix by $\tilde{\Gamma}_n:=(\tilde{\rho}_{ij})$
with $\tilde{\rho}_{ij}$ given as in (\ref{corr'}).

In this section we are interested in the limiting laws of the
coherence $L_n$ and $\tilde{L}_n$ of random matrix $X_n$, which are
defined to be the largest magnitude of the off-diagonal entries of
sample correlation matrices $\Gamma_n$ and $\tilde{\Gamma}_n$ respectively, see (\ref{mount}). The case of
$p \gg n$ is of particular interest to us. In such a setting,
some simulation studies about the distribution of $L_n$ were made in
Cai and Lv (2007), Fan and Lv (2008 and 2010). We now derive the
limiting laws of $L_n$ and $\tilde{L}_n$.

\newpage
We shall introduce another quantity that is useful for our technical
analysis. Define
\begin{eqnarray}\lbl{lidu}
J_n=\max_{1\leq i<j \leq p}{|(x_i-\mu)^T(x_j-\mu)|\over \sigma^2}.
\end{eqnarray}

We first state the law of large numbers
for $L_n$ for the case where the random entries $x_{ij}$ are bounded.

\begin{theorem}\lbl{army} Assume $|x_{11}|\leq C$ for a finite constant $C>0,$ and $p=p(n)\to \infty$ and $\log p=o(n)$ as $n\to\infty.$ Then
$\sqrt{n/\log p}\,L_n \to 2$ in probability as $n\to\infty.$
\end{theorem}

We now consider the case where $x_{ij}$ have finite exponential moments.
\begin{theorem}\lbl{morethan} Suppose $Ee^{t_0|x_{11}|^{\alpha}}< \infty$ for some $\alpha>0$ and $t_0>0.$ Set $\beta=\alpha/(4+\alpha).$  Assume  $p=p(n)\to \infty$ and $\log p=o(n^{\beta})$ as $n\to\infty.$ Then $\sqrt{n/\log p}\,L_n \to 2$ in probability as $n\to\infty.$
\end{theorem}

Comparing Theorems \ref{army} and \ref{morethan}, it can be seen that a stronger moment condition gives
a higher order of $p$ to make the law of large numbers for $L_n$
valid. Also, based on Theorem \ref{morethan}, if $Ee^{|x_{11}|^{\alpha}}< \infty$ for any $\alpha>0,$ then $\beta\to 1$, hence the order $o(n^{\beta})$ is close to $o(n),$ which is the order in Theorem \ref{army}.

We now consider the limiting distribution of $L_n$ after suitable
normalization.
\begin{theorem}\lbl{birthday} Suppose $Ee^{t_0|x_{11}|^{\alpha}}< \infty$ for some $0<\alpha\leq 2$ and $t_0>0.$ Set $\beta=\alpha/(4+\alpha).$  Assume  $p=p(n)\to \infty$ and $\log p=o(n^{\beta})$ as $n\to\infty.$ Then
$nL_n^2-4\log p + \log\log p$ converges weakly to an extreme
distribution of type I with
distribution function
\[
F(y)=e^{-{1\over \sqrt{8\pi}}e^{-y/2}}, \, y\in \mathbb{R}.
\]
\end{theorem}

\begin{remark} {\rm
Propositions \ref{father}, \ref{sister} and \ref{leave}
  show that the above three theorems are still valid if $L_n$ is
  replaced by either $\tilde{L}_n$ or $J_n/n,$ where $\tilde{L}_n$ is as in (\ref{mount}) and $J_n$ is as in (\ref{lidu}).
}
\end{remark}

In the case where $n$ and $p$ are comparable, i.e., $n/p\to \gamma \in
(0,\infty)$, Jiang (2004a) obtained the strong laws and
asymptotic distributions of the coherence $L_n$ of random
matrices. Several authors improved the results by sharpening the moment
assumptions, see, e.g., Li and Rosalsky (2006), Zhou (2007),  and Li,
Liu and Rosalsky (2009) where the same condition $n/p\to \gamma \in
(0,\infty)$ was imposed. Liu, Lin and Shao (2008)  showed that the
same results hold for $p\to\infty$ and $p=O(n^{\alpha})$ where
$\alpha$ is a constant.

In this paper, motivated by the applications mentioned earlier, we
are particularly interested in the case where both $n$ and $p$ are large and
$p=o(e^{n^{\beta}})$  while the entries of $X_n$ are i.i.d. with   a
certain moment condition. We also consider the case where the $n$ rows of $X_n$ form a random sample from $N_p(\mu, \Sigma)$ with  $\Sigma$ being a banded matrix. In particular, the entries of $X_n$ are not necessarily independent.    As shown in the above theorems and in Section \ref{dependent.sec} later, when $p\leq
e^{n^{\beta}}$ for a certain $\beta>0$, we obtain the strong laws and
limiting distributions of the coherence of random matrix
$X_n$. Presumably the results on high order $p=o(e^{n^{\beta}})$ need
stronger moment conditions than those for the case
$p=O(n^{\alpha})$. Ignoring the moment conditions, our results cover
those in  Liu, Lin and Shao (2008) as well as others aforementioned.

Theorem 1.2 in Jiang (2004a) states that if $n/p \to \gamma\in
(0, \infty)$  and  $E|\xi|^{30+\epsilon} <\infty$ for some $\epsilon
>0$, then for any $y \in \mathbb{R},$
\bea\lbl{preceding}
P\left(nL_n^2-4\log n+\log\log n \leq y\right) \to e^{-Ke^{-y/2}}
\eea
where $K=(\gamma^2\sqrt{8\pi})^{-1}$, as $n \to \infty$.
It is not difficult to see that Theorem \ref{birthday} implies Theorem 1.2 in Jiang (2004a) under condition that $n/p \to \gamma$ and $Ee^{t_0|x_{11}|^{\alpha}}< \infty$ for some $0<\alpha\leq 2$ and $t_0>0.$ In fact, write
\begin{eqnarray*}
& & nL_n^2-4\log n+\log\log n \\
& = & (nL_n^2-4\log p +\log\log p) + 4\log \frac{p}{n} +  \Big(\log\log n- \log\log p\Big).
\end{eqnarray*}
Theorem \ref{birthday} yields that $ nL_n^2-4\log p + \log\log p$
converges weakly to $F(y)=\exp^{-{1\over \sqrt{8\pi}}e^{-y/2}}$. Note
that since  $n/p \to \gamma$,
\[
4\log \frac{p}{n} \to -4\log \gamma \quad\mbox{and}\quad
\log(\log n)- \log \log p \to 0 .
\]
Now it follows from Slutsky's Theorem that $nL_n^2-4\log n +\log\log n$
converges weakly to $F(y + 4\log \gamma),$ which is exactly (\ref{preceding}) from Theorem
1.2 in Jiang (2004a).

\subsection{The Dependent Case}
\label{dependent.sec}

We now consider the case where the rows of random matrix $X_n$ are
drawn independently from a multivariate Gaussian distribution.
Let $X_n=(x_{ij})_{1\leq i \leq n, 1\leq j \leq p},$ where the $n$
rows are i.i.d. random vectors with distribution $N_p(\mu, \Sigma),$
where $\mu\in \mathbb{R}^p$ is  arbitrary in this section unless
otherwise specified. Let $(r_{ij})_{p\times p}$ be the correlation
matrix obtained from $\Sigma=(\sigma_{ij})_{p\times p}.$   As
mentioned in the introduction,  it is of interest to
test the hypothesis that the covariance matrix $\Sigma$ is banded, that is,
\begin{eqnarray}\lbl{Hypo}
H_0: \sigma_{ij}=0\ \mbox{for all } |i-j|\geq \tau
\end{eqnarray}
for a given integer $\tau\ge 1$. In order to construct a test, we study
in this section the asymptotic distributions of $L_{n,\tau}$ and
$\tilde{L}_{n,\tau}$ defined in (\ref{lag-k}) and (\ref{lag-k'}) respectively,
assuming the covariance matrix $\Sigma$ has desired banded structure
under the null hypothesis.
This case is much harder than the i.i.d. case considered in Section
\ref{iid.sec} because of the dependence.

For any $0<\delta<1,$ set
\begin{eqnarray}
\Gamma_{p,\delta}=\{1\leq i \leq p\,;\ |r_{ij}|>  1-\delta\ \mbox{for some}\  1\leq j \leq p\ \mbox{with}\ j\ne i\}.
\end{eqnarray}

\begin{theorem}\lbl{youngth}  Suppose, as $n\to\infty,$

(i) $p=p_n\to\infty$ with $\log p=o(n^{1/3})$;

(ii) $\tau=o(p^t)$ for any $t>0$;

(iii) for some $\delta\in (0,1)$, $|\Gamma_{p,\delta}|=o(p)$, which is particularly true if $\max_{1\leq i< j \leq p<\infty}|r_{ij}| \leq 1-\delta.$

\noindent Then, under $H_0$,
$nL_{n,\tau}^2-4\log p +\log\log p$ converges weakly to an extreme
distribution of type I with
distribution function
\[
F(y)=e^{-{1\over \sqrt{8\pi}}e^{-y/2}}, \, y\in \mathbb{R}.
\]
\end{theorem}
Similar to $J_n$ in (\ref{lidu}), we define
\begin{eqnarray}\lbl{take}
U_{n,\tau}=\max_{1\leq i<j \leq p,\, |i-j|\geq \tau}{|(x_i-\mu_i)^T(x_j-\mu_j)|\over \sigma_i\sigma_j}
\end{eqnarray}
where we write $x_i-\mu_i$
for $x_i-\mu_i e$ with $e=(1,1,\cdots, 1)^T \in \mathbb{R}^n$, $\mu=(\mu_1, \cdots,\mu_p)^T$ and $\sigma_i^2$'s are diagonal entries of $\Sigma.$

\begin{remark} {\rm
From Proposition \ref{correlated}, we know Theorem \ref{youngth} still holds if $L_{n,\tau}$ is replaced with $U_{n,\tau}$ defined  in (\ref{take}). In fact, by  the first paragraph in the proof of Theorem \ref{youngth}, to see if Theorem \ref{youngth} holds for $U_{n,\tau}$, we only need to consider the problem by assuming, w.l.o.g.,  $\mu=0$ and $\sigma_{i}$'s, the diagonal entries of $\Sigma$, are all equal to $1.$ Thus, by Proposition \ref{correlated}, Theorem \ref{youngth} holds when $L_{n,\tau}$ is replaced by $U_{n,\tau}.$
}
\end{remark}
Theorem \ref{youngth} implies immediately the following result.
\begin{coro} Suppose the conditions in Theorem \ref{youngth} hold, then $\sqrt{\frac{n}{\log p}}\,L_{n,\tau} \to 2$ in probability as $n\to\infty.$
\end{coro}
The assumptions (ii) and (iii) in Theorem \ref{youngth} are both essential. If one of them is violated, the conclusion may fail. The following two examples illustrate this point.

\begin{remark}{\rm
 Consider $\Sigma=I_{p}$ with $p=2n$ and $\tau=n.$ So conditions (i) and (iii) in Theorem \ref{youngth} hold, but (ii) does not. Observe
\begin{eqnarray*}
\Big\{(i,j);\, 1\leq i <  j \leq 2n,\, |i-j|\geq n\Big\}= n+(n-1)+\cdots + 1=\frac{n(n+1)}{2}\sim \frac{p^2}{8}
\end{eqnarray*}
as $n\to\infty.$ So $L_{n,\tau}$ is the maximum of roughly $p^2/8$  random variables, and the dependence of any two of such random variables are less than that appeared in $L_n$ in Theorem \ref{birthday}. The result in Theorem \ref{birthday} can be rewritten as
\begin{eqnarray*}
nL_n^2-2\log \frac{p^2}{2} + \log \log \frac{p^2}{2} -\log 8\ \ \mbox{converges weakly to}\ F
\end{eqnarray*}
as $n\to\infty.$ Recalling $L_n$ is the maximum of roughly $p^2/2$ weakly dependent random variables, replace $L_n$ with $L_{n,\tau}$ and $p^2/2$ with $p^2/8$ to have $nL_{n,\tau}^2-2\log \frac{p^2}{8} + \log \log \frac{p^2}{8} -\log 8$ converges weakly to $F,$ where $F$ is as in Theorem \ref{birthday}. That is,
\begin{eqnarray}\lbl{bulldozer}
(nL_{n,\tau}^2-4\log p + \log \log p) + \log 16\ \ \mbox{converges weakly to}\ F
\end{eqnarray}
as $n\to\infty$ (This can be done rigorously by following the proof of
Theorem \ref{birthday}). The difference between (\ref{bulldozer}) and
Theorem \ref{youngth} is evident.
}
\end{remark}

\begin{remark}{\rm
 Let $p=mn$ with integer $m\geq 2.$ We consider the $p\times p$ matrix $\Sigma=\ \mbox{diag}\, (H_n, \cdots, H_n)$ where there are $m$ $H_n$'s in the diagonal of $\Sigma$ and  all of the entries of the $n\times n$ matrix $H_n$ are equal to $1.$ Thus,  if $(\zeta_1, \cdots, \zeta_p)\sim N_p(0, \Sigma),$ then $\zeta_{ln+1}=\zeta_{ln+2}=\cdots=\zeta_{(l+1)n}$ for all $0\leq l \leq m-1$  and $\zeta_1, \zeta_{n+1}, \cdots, \zeta_{(m-1)n +1}$ are i.i.d. $N(0, 1)$-distributed random variables. Let $\{\zeta_{ij};\, 1\leq i \leq n, 1\leq j \leq m\}$ be i.i.d. $N(0, 1)$-distributed random variables. Then
\begin{eqnarray*}
(\underbrace{\zeta_{i1}, \cdots, \zeta_{i1}}_{n}, \underbrace{\zeta_{i2},\cdots, \zeta_{i2}}_{n}, \cdots, \underbrace{\zeta_{i\, m}, \cdots, \zeta_{i\, m}}_{n} )' \in \mathbb{R}^p,\ 1\leq i \leq n,
\end{eqnarray*}
 are i.i.d. random vectors with distribution $N_p(0, \Sigma).$ Denote the corresponding data matrix by $(x_{ij})_{n\times p}.$ Now, take $\tau=n$ and $m=[e^{n^{1/4}}].$ Notice $\Gamma_{p,\delta}=p$ for any $\delta>0.$ Since $p=mn,$ both (i) and (ii) in Theorem \ref{youngth} are satisfied, but (iii)  does not. Obviously,
\begin{eqnarray*}
L_{n,\tau}=\max_{1\leq i <  j \leq p,\, |i-j|\geq \tau}|\rho_{ij}|=\max_{1\leq i< j \leq m}|\hat{\rho}_{ij}|,
\end{eqnarray*}
where $\hat{\rho}_{ij}$ is obtained from $(\zeta_{ij})_{n\times m}$ as in (\ref{corr}) (note that the $mn$ entries of $(\zeta_{ij})_{n\times m}$ are i.i.d. with distribution $N(0,1)$). By Theorem \ref{birthday} on $\max_{1\leq i< j \leq m}|\hat{\rho}_{ij}|$, we have that $nL_{n,\tau}^2-4\log m+\log \log m$ converges weakly to $F,$ which is the same as the $F$  in Theorem \ref{youngth}. Set $\log_2x=\log \log x$ for $x> 1.$ Notice
\begin{eqnarray*}
nL_{n,\tau}^2-4\log m+\log_2 m & = & nL_{n,\tau}^2-4\log p + 4\log n + \log _2m\\
& \sim & (nL_{n,\tau}^2-4\log p+\log_2 p) + 4\log n
\end{eqnarray*}
since $p=mn$ and $\log_2p - \log_2m \to 0.$ Further, it is easy to check that $4\log n-16\log_2p \to 0$. Therefore, the previous  conclusion is equivalent to that
\begin{eqnarray}\lbl{ahaa}
(nL_{n,\tau}^2-4\log p+\log \log p) + 16\log \log p \ \mbox{converges weakly to}\ F
\end{eqnarray}
as $n \to \infty.$ This is different from the conclusion of Theorem \ref{youngth}.
}
\end{remark}

\section{Testing the Covariance Structure}
\label{testing.sec}

The limiting laws derived in the last section have immediate
statistical applications. Testing the covariance structure of a high
dimensional random variable is an important
problem in statistical inference. In particular,
as aforementioned, in econometrics when testing certain economic
theories and in time series analysis in general
 it is of significant interest to test the hypothesis that the
 covariance matrix $\Sigma$ is banded. That is, the variables have nonzero
correlations only up to a certain lag $\tau$.
The limiting distribution of $L_{n,\tau}$ obtained in Section
\ref{limit.sec} can be readily used to construct a test for the
bandedness of the covariance matrix in the Gaussian case.

Suppose we observe independent and identically distributed $p$-variate
Gaussian variables $\mathbf{Y}_{1},\ldots ,\mathbf{Y}_{n}$ with mean
$\mu_{p\times 1}$, covariance matrix $\Sigma _{p\times p}=(\sigma_{ij})$ and
correlation matrix $R _{p\times p}=(r_{ij})$. For a given integer
$\tau\ge 1$ and a given significant level $0<\alpha < 1$,  , we wish
to test the hypotheses
\be
\label{hypothesis}
H_0: \; \sigma_{i,j} = 0 \;\;\mbox{for all $|i-j|\ge \tau$}
\;\;\mbox{versus} \;\;
H_a: \sigma_{i,j} \neq 0 \;\;\mbox{for some $|i-j|\ge \tau$}.
\ee
A case of special interest is $\tau =1$, which corresponds to testing
independence of the Gaussian random variables.
The asymptotic distribution of  $L_{n,\tau}$ derived in Section
\ref{dependent.sec} can be used to construct a convenient test
statistic for testing the hypotheses in (\ref{hypothesis}).

Based on the asymptotic result given in Theorem \ref{youngth} that
\bea\lbl{sparrow}
P\left(n L_{n,\tau}^2 - 4\log p + \log\log p \le y\right) \goto e^{-{1\over
    \sqrt{8\pi}} e^{-y/2}},
\eea
we define a test for testing the hypotheses in (\ref{hypothesis}) by
\be
\label{test}
T = I\Big(L_{n,\tau}^2 \ge n^{-1} (4\log p - \log\log p - \log (8\pi) - 2 \log\log(1-\alpha)^{-1})\Big).
\ee
That is, we reject the null hypothesis $H_0$ whenever
\[
L_{n,\tau}^2 \ge n^{-1} \Big(4\log p - \log\log p - \log (8\pi) - 2 \log\log(1-\alpha)^{-1}\Big).
\]
Note that for $\tau=1$, $L_{n,\tau}$ reduces to $L_n$ and the test is then
based on the coherence $L_n$.

\begin{theorem}
Under the conditions of Theorem \ref{youngth},
the test $T$ defined in (\ref{test}) has size $\alpha$ asymptotically.
\end{theorem}
This result is a direct consequence of (\ref{sparrow}).
\begin{remark} {\rm
For testing independence, another natural approach is to build a test
based on the largest eigenvalue $\lambda_{\max}$ of the sample
correlation matrix. However, the limiting distribution of the largest
eigenvalue $\lambda_{\max}$ is unknown even for the case $p/n\to c,$ a finite and positive constant. For
$\tau\ge 2$, the eigenvalues are not useful for testing bandedness of the
covariance matrix.
}
\end{remark}

\section{Construction of Compressed Sensing Matrices}
\label{CS.sec}

As mentioned in the introduction, an important problem in compressed
sensing is the construction of measurement matrices $X_{n\times p}$
which enables the precise recovery of a sparse signal $\beta$ from
linear measurements $y = X\beta$  using an efficient recovery
algorithm. Such a measurement matrix $X$ is difficult to
construct deterministically. It has been shown that randomly generated
matrix $X$ can satisfy the so called RIP condition with high
probability.

The best known example is perhaps $n\times p$ random matrix
$X$ whose entries $x_{i,j}$ are iid normal variables
\be
\label{normal.CS}
x_{i,j}\stackrel{iid}{\sim} N(0, n^{-1}).
\ee
Other examples include generating $X=(x_{i,j})$ by Bernoulli
random variables
\be
\label{bernoulli.CS}
x_{i,j}=\left\{
\begin{array}{ll}
1/\sqrt{n} & \quad \mbox{with probability $\hf$;}\\
-1/\sqrt{n} &\quad \mbox{with probability $\hf$}
\end{array}
\right.
\ee
or more sparsely by
\be
\label{sparse.CS}
x_{i,j}=\left\{
\begin{array}{ll}
\sqrt{3/n} & \quad \mbox{with probability $1/6$;}\\
0 & \quad \mbox{with probability $2/3$;}\\
-\sqrt{3/n} &\quad \mbox{with probability $1/6$.}
\end{array}
\right.
\ee
These random matrices are shown to satisfy the RIP conditions with
high probability. See Achlioptas (2001) and Baraniuk, et al. (2008).

In addition to RIP, another commonly used condition is the mutual
incoherence property (MIP) which
requires the pairwise correlations among the column vectors of $X$ to
be small. In compressed sensing $\tilde{L}_n$ (instead of $L_n$) is commonly
used.
It has been shown that the condition
\be
\label{sharp.condition1}
(2k-1)\tilde{L}_n < 1
\ee
ensures the exact recovery of $k$-sparse signal $\beta$ in the
noiseless case where $y=X\beta$, and stable recovery of sparse signal in the
noisy case where
\[
y=X\beta + z.
\]
Here $z$ is an error vector, not necessarily random.
The MIP (\ref{sharp.condition1}) is  a very
desirable property. When  the measurement matrix $X$ satisfies
(\ref{sharp.condition1}), the constrained $\ell_1$ minimizer can be shown to be
exact in the noiseless case and near-optimal in the noisy case.
Under the MIP condition, the analysis of $\ell_1$ minimization methods
is also  particularly simple. See, e.g., Cai, Wang and Xu (2010b).

The results given in Theorems \ref{army} and \ref{morethan} can be
used to show how likely a random matrix satisfies the MIP condition
(\ref{sharp.condition1}). Under the conditions of either Theorem
\ref{army} or Theorem  \ref{morethan},
\[
\tilde{L}_n \sim 2\sqrt{\log p \over n}.
\]
So in order for the MIP condition (\ref{sharp.condition1}) to hold,
roughly the sparsity $k$ should satisfy
\[
k < {1\over 4} \sqrt{n\over \log p}.
\]
In fact we have the following more precise result which is proved in
Section \ref{proof.sec}.

\begin{prop}\lbl{applesensing} Let $X_n=(x_{ij})_{n\times p}$  where $x_{ij}$'s
 are i.i.d. random variables with mean $\mu$, variance $\sigma^2>0$ and  $Ee^{t_0|x_{11}|^2}< \infty$ for some $t_0>0.$ Let $\tilde{L}_n$ be as in (\ref{mount}). Then $P(\tilde{L}_n \geq t) \leq 3 p^2 e^{-n g(t)}$
where $g(t) =\min\{I_1(t/2),\, I_2(1/2)\}> 0$ for any $t>0$  and
\beaa
I_1(x)=\sup_{\theta\in \mathbb{R}}\{\theta x-\log Ee^{\theta \xi \eta}\}\ \mbox{and}\ I_2(x)=\sup_{\theta \in \mathbb{R}}\{\theta x-\log Ee^{\theta \xi^2}\}.
\eeaa
and $\xi, \eta, (x_{11}-\mu)/\sigma$ are i.i.d.
\end{prop}

We now consider the three particular random matrices mentioned in the
beginning of this section.

\noindent\textbf{Example 1}. Let $x_{11}\sim N(0, n^{-1})$ as in (\ref{normal.CS}). In this case, according to the above proposition, we have
\bea\lbl{chop}
P\left((2k-1)\tilde{L}_n < 1\right) \geq 1- 3 p^2\exp\Big\{-\frac{n}{12(2k-1)^2}\Big\}
\eea
for all $n\geq 2$ and $k\geq 1.$ The verification of this example together with the next two are given in the Appendix.

\noindent\textbf{Example 2}. Let $x_{11}$ be such that $P(x_{11}=\pm 1/\sqrt{n})=1/2$ as in (\ref{bernoulli.CS}). In this case, we have
\bea\lbl{karaok}
P\left((2k-1)\tilde{L}_n < 1\right) \geq 1- 3 p^2\exp\Big\{-\frac{n}{12(2k-1)^2}\Big\}
\eea
for all $n\geq 2$ and $k\geq 1.$

\noindent\textbf{Example 3}. Let $x_{11}$ be such that $P(x_{11}=\pm \sqrt{3/n})=1/6$ and $P(x_{11}=0)=2/3$ as in (\ref{sparse.CS}). Then
\bea\lbl{roasted}
P\left((2k-1)\tilde{L}_n < 1\right) \geq 1- 3 p^2\exp\Big\{-\frac{n}{12(2k-1)^2}\Big\}
\eea
for all $n\geq 2$ and $k\geq 2.$

\begin{remark}{\rm One can see from the above that (\ref{chop}) is true for all of the three examples with different restrictions on $k.$ In fact this is always the case as long as $Ee^{t_0|x_{11}|^2}< \infty$ for some $t_0>0,$ which can be seen from Lemma \ref{solar.system}.
}
\end{remark}

\begin{remark}{\rm
Here we would like to point out an error on pp. 801 of Donoho (2006b)
and pp. 2147 of Candes and Plan (2009) that the coherence of a random
matrix with i.i.d. Gaussian entries is about $2\sqrt{\log p \over n}$,
not $\sqrt{2\log p \over n}$.
}
\end{remark}






\section{Discussion and Comparison with Related Results}\lbl{discussion.sec}

This paper studies the limiting laws of the largest magnitude of the
off-diagonal entries of the sample correlation matrix in the
high-dimensional setting. Entries of other types of random
matrices have been studied in the literature, see, e.g.,
Diaconis, Eaton and Lauritzen (1992), and Jiang (2004a, 2005, 2006,
2009).
Asymptotic properties of the eigenvalues of the sample correlation
matrix  have also been studied when both $p$ and $n$ are large
and proportional to each other. For instance, it is proved in Jiang (2004b)
 that the
empirical distributions of the eigenvalues of the sample correlation
matrices converge  to the Marchenko-Pastur law; the largest and
smallest eigenvalues satisfy certain law of large numbers.
However, the high-dimensional case of $p\gg n$ remains an open problem.

The motivations of our current
work consist of the applications to testing covariance structure  and
construction of compressed sensing matrices in the ultra-high
dimensional setting where the dimension $p$ can be as large as
$e^{n^\beta}$ for some $0<\beta<1$. The setting is different from
those considered in the earlier literature such as Jiang (2004), Zhou
(2007), Liu, Lin and Shao (2008), and Li, Liu and Rosalsky (2009).
Our main theorems and techniques
are different from those mentioned above in the following two aspects:

\begin{itemize}
\item[(a)] Given $n\to\infty,$ we push the size of $p$ as large as we
  can to make the law of large numbers and limiting results on $L_n$
  and $\tilde{L}_n$ valid. Our current theorems say that, under some moment conditions, these results hold as long as $\log p =o(n^{\beta})$ for a certain $\beta>0$.

\item[(b)] We study $L_n$ and $\tilde{L}_n$ when the $p$ coordinates
  of underlying multivariate distribution are not i.i.d. Instead, the
  $p$ coordinates follow a multivariate normal distribution $N_p(\mu,
  \Sigma)$ with $\Sigma$ being banded and $\mu$ arbitrary. Obviously, the $p$ coordinates are dependent. The proofs of our theorems are more subtle and involved than those in the earlier papers. In fact, we have to consider the dependence structure of $\Sigma$ in detail, which is more complicated than the independent case. See Lemmas \ref{webster}, \ref{yahoo} and \ref{yahoo1}.
\end{itemize}

Liu, Lin and Shao (2008) introduced a statistic for testing
independence that is different from
$L_n$ and $\tilde{L}_n$ to improve the convergence speed of the two
statistics under the constraint $c_1 n^\alpha \le p \le c_2 n^\alpha$
for some constants $c_1, c_2, \alpha > 0$. In this paper, while
pushing the order of $p$ as large as possible to have the limit
theorems, we focus on the behavior of $L_n$ and $\tilde{L}_n$
only. This is because $L_n$ and $\tilde{L}_n$ are specifically
used in some applications such as compressed sensing. On the other hand, we also consider a more
general testing problem where one wishes to test the bandedness of the
covariance matrix $\Sigma$ in $N_p(\mu, \Sigma)$ while allowing $\mu$
to be arbitrary. We propose the statistic $L_{n, \tau}$ in (\ref{lag-k})
and derive its law of large numbers and its limiting distribution. To
our knowledge, this is new in the literature. It is
interesting to explore the possibility of improving the convergence
speed by modifying $L_{n, \tau}$ as that of $L_n$ in Liu, Lin and Shao
(2008). We leave this as future work.

\section{Proofs}
\label{proof.sec}

In this section we prove Theorems \ref{army} - \ref{youngth}. The letter $C$ stands for a constant and may vary from place to place throughout this section. Also, we sometimes write $p$ for $p_n$ if there is no confusion.
For any square matrix $A=(a_{i,j}),$ define $\mn A\mn=\max_{1\leq i\ne j \leq n}|a_{i,j}|;$ that is, the maximum of the absolute values of the off-diagonal entries of $A.$

We begin by collecting a few essential technical lemmas in Section
\ref{Technical} without proof. Other technical lemmas used in the
proofs of the main results are proved in the Appendix.

\subsection{Technical Tools}\lbl{Technical}

\begin{lemma}\lbl{chicago}(Lemma 2.2 from Jiang (2004a)) Recall $x_i$ and $\Gamma_n$ in (\ref{corr}). Let $h_i=\|x_i-\bar{x}_i\|/\sqrt{n}$ for each $i.$ Then
\beaa\lbl{newadd}
\mn n\Gamma_n-X_n^TX_n\mn \leq (b_{n,1}^2+2b_{n,1})W_nb_{n,3}^{-2}+nb_{n,3}^{-2}b_{n,4}^2,
\eeaa
where
\beaa
& & b_{n,1}=\max_{1\leq i \leq p}|h_i -1|,\ \ \ W_n=\max_{1\leq i < j \leq p}|x_i^Tx_j|,\ \ \ b_{n,3}=\min_{1\leq i \leq p}h_i, \ \ \ b_{n,4}=\max_{1\leq i \leq p}|\bar{x}_i|.
\eeaa
\end{lemma}

The following Poisson approximation result
is essentially a special case of
Theorem $1$ from Arratia et al. (1989).
\begin{lemma}\label{stein} Let $I$ be an index set and $\{B_{\alpha}, \alpha\in I\}$ be a set of subsets of $I,$ that is, $B_{\alpha}\subset I$ for each $\alpha \in I.$  Let also $\{\eta_{\alpha}, \alpha\in I\}$ be random variables. For a given $t\in \mathbb{R},$ set $\lambda=\sum_{\alpha\in I}P(\eta_{\alpha}>t).$ Then
\beaa
|P(\max_{\alpha \in I}\eta_{\alpha} \leq t)-e^{-\lambda}| \leq (1\wedge \lambda^{-1})(b_1+b_2+b_3)
\eeaa
where
\beaa
& & b_1=\sum_{\alpha \in I}\sum_{\beta \in B_{\alpha}}P(\eta_{\alpha} >t)P(\eta_{\beta} >t),\\
& & b_2=\sum_{\alpha \in I}\sum_{\alpha\ne \beta \in B_{\alpha}}P(\eta_{\alpha} >t, \eta_{\beta} >t),\\
 & & b_3=\sum_{\alpha \in I}E|P(\eta_{\alpha} >t|\sigma(\eta_{\beta}, \beta \notin B_{\alpha})) - P(\eta_{\alpha} >t)|,
\eeaa
and $\sigma(\eta_{\beta}, \beta \notin B_{\alpha})$ is the $\sigma$-algebra generated by $\{\eta_{\beta}, \beta \notin B_{\alpha}\}.$
In particular, if $\eta_{\alpha}$ is independent of $\{\eta_{\beta}, \beta \notin B_{\alpha}\}$ for each $\alpha,$ then $b_3=0.$
\end{lemma}

The following conclusion is Example 1 from Sakhanenko (1991). See also Lemma 6.2 from Liu et al (2008).
\begin{lemma}\lbl{shao} Let $\xi_i, 1\leq i \leq n,$ be independent random variables with $E\xi_i=0.$ Put
\begin{eqnarray*}
s_n^2=\sum_{i=1}^nE\xi_i^2,\ \ \ \varrho_n=\sum_{i=1}^nE|\xi_i|^3,\ \ \ S_n=\sum_{i=1}^n\xi_i.
\end{eqnarray*}
Assume $\max_{1\leq i \leq n}|\xi_i| \leq c_ns_n$ for some $0<c_n \leq 1.$ Then
\begin{eqnarray*}
P(S_n\geq xs_n) =e^{\gamma(x/s_n)}(1-\Phi(x))(1+\theta_{n,x}(1+x)s_n^{-3}\varrho_n)
\end{eqnarray*}
for $0<x \leq 1/(18c_n),$ where $|\gamma(x)| \leq 2x^3\varrho_n$ and $|\theta_{n,x}|\leq 36.$
\end{lemma}
The following are moderate deviation results from Chen (1990), see also Chen (1991), Dembo and Zeitouni (1998) and Ledoux (1992). They are a special type of large deviations.
\begin{lemma}\lbl{Xia} Suppose $\xi_1, \xi_2, \cdots$ are i.i.d. r.v.'s with  $E\xi_1=0$ and $E\xi_1^2=1.$ Put $S_n=\sum_{i=1}^n\xi_i.$

(i) Let $0< \alpha \leq 1$ and $\{a_n;\, n\geq 1\}$ satisfy that  $a_n\to+\infty$ and $a_n=o\big(n^{\frac{\alpha}{2(2-\alpha)}}\big).$ If $Ee^{t_0|\xi_1|^{\alpha}}<\infty$ for some $t_0>0,$ then
\begin{eqnarray}\lbl{reading}
\lim_{n\to\infty}\frac{1}{a_n^2}\log P\Big(\frac{S_n}{\sqrt{n}a_n}\geq u\Big)=-\frac{u^2}{2}
\end{eqnarray}
for any $u>0.$

(ii) Let $0<\alpha<1$ and $\{a_n;\, n\geq 1\}$ satisfy that  $a_n\to+\infty$ and $a_n=O\big(n^{\frac{\alpha}{2(2-\alpha)}}\big).$ If $Ee^{t|\xi_1|^{\alpha}}<\infty$ for all $t>0,$ then (\ref{reading}) also holds.
\end{lemma}

\subsection{Proofs of Theorems \ref{army} and \ref{morethan}}\label{LLN}

Recall that a sequence of random variables $\{X_n;\, n\geq 1\}$ are said to be {\it tight} if, for any $\epsilon>0,$ there is a constant $K>0$ such that $\sup_{n\geq 1}P(|X_n|\geq K)< \epsilon.$ Obviously, $\{X_n;\, n\geq 1\}$ are tight if for some $K>0$, $\lim_{n\to\infty}P(|X_n|\geq K) \to 0.$ It is easy to check that
\begin{eqnarray}\lbl{tight}
& & \mbox{if}\  \{X_n;\, n\geq 1\}\ \mbox{are tight, then for any sequence of constants}\ \{\epsilon_n;\ n\geq 1\}\ \ \ \ \ \ \ \ \ \ \ \ \ \ \ \ \ \nonumber\\
& & \mbox{with}\ \lim_{n\to\infty}\epsilon_n = 0,\    \mbox{we have}\ \epsilon_nX_n \to 0 \ \mbox{in probability as} \ n\to\infty.
\end{eqnarray}

Reviewing the notation $b_{n,i}$'s defined in Lemma \ref{chicago}, we have the following properties.
\begin{lemma}\lbl{palo} Let $\{x_{ij};\, i\geq 1, j\geq 1\}$ be i.i.d. random variables with  $Ex_{11}=0$ and $Ex_{11}^2=1.$  Then, $b_{n,3} \to 1$ in probability as $n\to \infty,$ and
$\{\sqrt{n/\log p}\, b_{n,1}\}$ and $\{\sqrt{n/\log p}\, b_{n,4}\}$ are  tight provided one of the following conditions holds:

(i) $|x_{11}|\leq C$ for some constant $C>0$, $p_n\to\infty$ and $\log p_n =o(n)$ as $n\to\infty$;

(ii) $Ee^{t_0|x_{11}|^{\alpha}}< \infty$ for some $0<\alpha \leq 2$ and $t_0>0,$ and $p_n\to\infty$ and $\log p_n=o(n^{\beta})$ as $n\to\infty,$ where $\beta=\alpha/(4-\alpha).$
\end{lemma}

\begin{lemma}\lbl{sweet} Let $\{x_{ij};\, i\geq 1, j\geq 1\}$ be i.i.d. random variables with $|x_{11}|\leq C$ for a finite constant $C>0$, $Ex_{11}=0$ and $E(x_{11}^2)=1.$  Assume $p=p(n)\to \infty$ and $\log p=o(n)$ as $n\to\infty.$  Then, for any $\epsilon>0$ and a sequence of positive numbers $\{t_n\}$ with limit $t>0,$
\begin{eqnarray*}
\Psi_n:= E\Big\{ P^1\Big(|\sum_{k=1}^nx_{k1}x_{k2}| > t_n\sqrt{n\log p}\,\Big)^2\Big\}=O\left(\frac{1}{p^{t^2-\epsilon}}\right)
\end{eqnarray*}
as $n\to\infty$, where $P^1$ stands for the conditional probability given $\{x_{k1},\, 1\leq k \leq n\}.$
\end{lemma}

\begin{lemma}\lbl{sweet1} Suppose $\{x_{ij};\, i\geq 1, j\geq 1\}$ are i.i.d. random variables with  $Ex_{11}=0,\,  E(x_{11}^2)=1$ and $Ee^{t_0|x_{11}|^{\alpha}}< \infty$ for some $t_0>0$ and $\alpha>0.$ Assume $p=p(n)\to \infty$ and $\log p=o(n^{\beta})$ as $n\to\infty,$ where $\beta=\alpha/(4+\alpha).$  Then, for any $\epsilon>0$ and a sequence of positive numbers $\{t_n\}$ with limit $t>0,$
\begin{eqnarray*}
\Psi_n:= E\Big\{ P^1\Big(|\sum_{k=1}^nx_{k1}x_{k2}| > t_n\sqrt{n\log p}\,\Big)^2\Big\}=O\left(\frac{1}{p^{t^2-\epsilon}}\right)
\end{eqnarray*}
as $n\to\infty$, where $P^1$ stands for the conditional probability given $\{x_{k1},\, 1\leq k \leq n\}.$
\end{lemma}
Lemmas \ref{palo}, \ref{sweet}, and \ref{sweet1} are proved in the Appendix.

\begin{prop}\lbl{father} Suppose the conditions in Lemma \ref{sweet} hold with $X_n=(x_{ij})_{n\times p}=(x_1, \cdots, x_p).$ Define $W_n=\max_{1\leq i<j \leq p}|x_i^Tx_j|=\max_{1\leq i<j \leq p}\left|\sum_{k=1}^nx_{ki}x_{kj}\right|.$  Then
\begin{eqnarray*}
\frac{W_n}{\sqrt{n\log p}} \to  2
\end{eqnarray*}
in probability as $n\to\infty.$
\end{prop}
\noindent\textbf{Proof}. We first prove
\begin{eqnarray}\lbl{brook}
\lim_{n\to\infty}P\Big(\frac{W_n}{\sqrt{n\log p}} \geq  2+2\epsilon\Big)= 0
\end{eqnarray}
for any $\epsilon>0.$ First, since $\{x_{ij};\, i\geq 1, j\geq 1\}$ are i.i.d., we have
\begin{eqnarray}\lbl{grass}
P(W_n \geq (2+2\epsilon)\sqrt{n\log p}) \leq \binom{p}{2}\cdot P\Big(\Big|\sum_{k=1}^nx_{k1}x_{k2}\Big|\geq (2+2\epsilon)\sqrt{n\log p}\Big)
\end{eqnarray}
for any $\epsilon>0.$ Notice $E(|x_{11}x_{12}|^2)=E(|x_{11}|^2)\cdot E(|x_{12}|^2)=1.$ By (i) of Lemma \ref{Xia}, using conditions  $Ee^{|x_{11}x_{12}|}<\infty$ and $\log p=o(n)$ as $n\to\infty$, we obtain
\begin{eqnarray}\lbl{sky}
P\left(\left|\sum_{k=1}^nx_{k1}x_{k2}\right|\geq (2+2\epsilon)\sqrt{n\log p}\right) \leq \exp\left(-\frac{(2+\epsilon)^2}{2}\log p\right)\leq \frac{1}{p^{2+\epsilon}}
\end{eqnarray}
as $n$ is sufficiently large. The above two assertions conclude
\begin{eqnarray}\lbl{kick}
P(W_n \geq (2+2\epsilon)\sqrt{n\log p}) \leq \frac{1}{p^{\epsilon}}\to 0
\end{eqnarray}
as $n\to \infty.$ Thus (\ref{brook}) holds. Now, to finish the proof, we only need to show
\begin{eqnarray}\lbl{tree}
\lim_{n\to\infty}P\Big(\frac{W_n}{\sqrt{n\log p}} \leq  2-\epsilon\Big)= 0
\end{eqnarray}
for any $\epsilon>0$ small enough.

Set $a_n=(2-\epsilon)\sqrt{n\log p}$ for $0<\epsilon <2$ and
\begin{eqnarray*}
y_{ij}^{(n)}=\sum_{k=1}^nx_{ki}x_{kj}
\end{eqnarray*}
for $1\leq i, j \leq n.$ Then $W_n=\max_{1\leq i < j \leq p}|y_{ij}^{(n)}|$ for all $n\geq 1.$

Take $I=\{(i,j);\ 1\leq i< j \leq p\}.$ For $u =(i,j) \in I,$ set $B_{u}=\{(k,l)\in I;\ \mbox{one of}\ k\ \mbox{and}\ l =i\ \mbox{or}\ j,\ \mbox{but}\ (k,l)\ne u\},\ \eta_{u}=|y_{ij}^{(n)}|,\ t=a_n$ and $A_{u}=A_{ij}=\{|y_{ij}^{(n)}| > a_n\}.$ By the i.i.d. assumption on $\{x_{ij}\}$ and Lemma \ref{stein},
\bea\lbl{season}
P(W_n\leq a_n) \leq e^{-\lambda_n} +b_{1,n}+b_{2,n}
\eea
where
\bea\lbl{sheep}
\lambda_n=\frac{p(p-1)}{2}P(A_{12}),\ b_{1,n}\leq 2p^3P(A_{12})^2\ \mbox{and}\ b_{2,n} \leq 2p^3P(A_{12}A_{13}).
\eea
Remember that $y_{12}^{(n)}$ is a sum of i.i.d. bounded random variables with mean $0$ and variance $1.$ By (i) of Lemma \ref{Xia}, using conditions  $Ee^{t|x_{11}x_{12}|}<\infty$ for any $t>0$ and $\log p=o(n)$ as $n\to\infty$, we know
\begin{eqnarray}\lbl{hurdle}
\lim_{n\to\infty}\frac{1}{\log p}\log P(A_{12}) = -\frac{(2-\epsilon)^2}{2}
\end{eqnarray}
for any $\epsilon \in (0, 2).$ Noticing $2-2\epsilon< (2-\epsilon)^2/2 < 2-\epsilon$ for $\epsilon \in (0, 1),$ we have that
\begin{eqnarray}\lbl{Yin}
\frac{1}{p^{2-\epsilon}}\leq P(A_{12})\leq \frac{1}{p^{2-2\epsilon}}
\end{eqnarray}
as $n$ is sufficiently large. This implies
\begin{eqnarray}\lbl{laughing}
e^{-\lambda_n} \leq e^{-p^{\epsilon}/3}\ \ \mbox{and}\ \ \ b_{1,n} \leq \frac{2}{p^{1-4\epsilon}}
\end{eqnarray}
for $\epsilon \in (0, 1/4)$ as $n$ is large enough. On the other hand, by independence
\begin{eqnarray}
P(A_{12}A_{13})
& = & P(|y_{12}^{(n)}| > a_n, |y_{13}^{(n)}| > a_n)\lbl{sleeping}\\
& = & E\big\{ P^1(|\sum_{k=1}^nx_{k1}x_{k2}| > a_n)^2\big\}\nonumber
\end{eqnarray}
where $P^1$ stands for the conditional probability given $\{x_{k1},\, 1\leq k \leq n\}.$ By Lemma \ref{sweet},
\begin{eqnarray}\lbl{Nike}
P(A_{12}A_{13})
 \leq   \frac{1}{p^{4-4\epsilon}}
\end{eqnarray}
for any $\epsilon>0$ as $n$ is sufficiently large. Therefore, taking $\epsilon \in (0,1/4),$ we have
\begin{eqnarray}\lbl{littleshow}
b_{2,n} \leq 2p^3P(A_{12}A_{13}) \leq \frac{2}{p^{1-4\epsilon}} \to 0
\end{eqnarray}
as $n\to\infty.$ This together with (\ref{season}) and (\ref{laughing}) concludes (\ref{tree}).\ \ \ \ \ \ \ \ $\blacksquare$\\

\begin{prop}\lbl{sister} Suppose the conditions in Lemma \ref{sweet1} hold. Let $W_n$ be as in Lemma \ref{father}.  Then
\begin{eqnarray*}
\frac{W_n}{\sqrt{n\log p}} \to  2
\end{eqnarray*}
in probability as $n\to\infty.$
\end{prop}
The proof of Proposition \ref{sister} is similar to that of
Proposition \ref{father}. Details are given in the Appendix.

\bigskip
\noindent\textbf{Proof of Theorem \ref{army}}. First, for constants $\mu_i \in \mathbb{R}$ and $\sigma_i > 0,\ i=1,2,\cdots,p,$ it is easy to see that matrix $X_n=(x_{ij})_{n\times p}=(x_1, x_2, \cdots, x_p)$ and $(\sigma_1x_1+\mu_1e, \sigma_2 x_2+\mu_2 e,\cdots, \sigma_p x_p+\mu_p e)$ generate the same sample correlation matrix $\Gamma_n=(\rho_{ij})$, where  $\rho_{ij}$ is as in (\ref{corr}) and $e=(1,\cdots, 1)' \in \mathbb{R}^n$. Thus, w.l.o.g., we prove the theorem next by assuming that $\{x_{ij};\, 1\leq i \leq n, 1\leq j \leq p\}$ are i.i.d. random variables with mean zero and variance $1.$

By Proposition \ref{father}, under condition $\log p=o(n),$
\begin{eqnarray}\lbl{five1}
\frac{W_n}{\sqrt{n\log p}} \to 2
\end{eqnarray}
in probability as $n\to\infty.$ Thus, to prove the theorem, it is enough to show
\begin{eqnarray}\lbl{preschool1}
\frac{nL_n-W_n}{\sqrt{n\log p}} \to 0
\end{eqnarray}
in probability as $n\to\infty.$   From Lemma \ref{chicago},
\begin{eqnarray*}
|nL_n-W_n| \leq \mn n\Gamma_n-X_n^TX_n\mn \leq (b_{n,1}^2+2b_{n,1})W_nb_{n,3}^{-2}+nb_{n,3}^{-2}b_{n,4}^2.
\end{eqnarray*}
By (i) of Lemma \ref{palo}, $b_{n,3} \to 1$ in probability as $n\to\infty,\ \{\sqrt{n/\log p}\, b_{n,1}\}$ and $\{\sqrt{n/\log p}\, b_{n,4}\}$ are all tight. Set $b_{n,1}'=\sqrt{n/\log p}\, b_{n,1}$ and $b_{n,4}'=\sqrt{n/\log p}\, b_{n,4}$ for all $n\geq 1.$ Then $\{b_{n,1}'\}$ and $\{b_{n,4}'\}$ are both tight. It follows that
\begin{eqnarray*}
\frac{|nL_n-W_n|}{\sqrt{n\log p}}\leq \sqrt{\frac{\log p}{n}}\,
\left(\sqrt{\frac{\log p}{n}}b_{n,1}'^2 + 2b_{n,1}'\right)\cdot\frac{W_n}{\sqrt{n\log p}}\cdot b_{n,3}^{-2} + \sqrt{\frac{\log p}{n}}\, b_{n,3}^{-2}b_{n,4}'^2,
\end{eqnarray*}
 which concludes (\ref{preschool1}) by (\ref{tight}).\ \ \ \ \ \ \ $\blacksquare$\\

\noindent\textbf{Proof of Theorem \ref{morethan}}. In the proof of Theorem \ref{army}, replace ``Proposition \ref{father}" with ``Proposition \ref{sister}" and ``(i) of Lemma \ref{palo}" with ``(ii) of Lemma \ref{palo}", keep all other statements the same, we then get the desired result.\ \ \ \ \ \ \ $\blacksquare$\\

\noindent\textbf{Proof of Proposition \ref{applesensing}}. Recall the definition of  $\tilde{L}_n$ in (\ref{mount}), to prove the conclusion, w.l.o.g., we assume $\mu=0$ and $\sigma^2=1.$ Evidently, by the i.i.d. assumption,
\begin{eqnarray}\lbl{insult}
P(\tilde{L}_n \geq t) & \leq & \frac{p^2}{2}P\Big(\frac{|x_1'x_2|}{\|x_1\|\cdot \|x_2\|} \geq t\Big)\nonumber\\
 & \leq & \frac{p^2}{2}P\Big(\frac{|x_1'x_2|}{n} \geq \frac{t}{2}\Big) + \frac{p^2}{2}\cdot 2 P\left(\frac{\|x_1\|^2}{n}\leq \frac{1}{2} \right)
\end{eqnarray}
where the event $\{\|x_{11}\|^2/n>1/2,\, \|x_{12}\|^2/n>1/2\}$ and its complement are used to get the last inequality. Since $\{x_{ij};\, i\geq 1,\, j\geq 1\}$ are i.i.d., the condition $Ee^{t_0|x_{11}|^2}< \infty$ implies $Ee^{t_0'|x_{11}x_{12}|}< \infty$ for some $t_0'>0.$ By the Chernoff bound (see, e.g., p. 27 from Dembo and Zeitouni (1998)) and noting that $E(x_{11}x_{12})=0$ and $E x_{11}^2=1,$ we have
\beaa
P\Big(\frac{|x_1'x_2|}{n} \geq \frac{t}{2}\Big)\leq 2 e^{-nI_1(t/2)}\ \mbox{and}\ P\left(\frac{\|x_1\|^2}{n}\leq \frac{1}{2} \right)\leq 2e^{-nI_2(1/2)}
\eeaa
for any $n\geq 1$ and $t>0,$ where the following facts about rate functions $I_1(x)$ and $I_2(y)$ are used:

(i) $I_1(x)=0$ if and only if $x=0$; $I_2(y)=0$ if and only if $y=1;$

(ii) $I_1(x)$ is non-decreasing on $A:=[0, \infty)$ and non-increasing on $A^c$. This is also true for  $I_2(y)$ with $A=[1, \infty).$

\noindent These and (\ref{insult}) conclude
\beaa
P(\tilde{L}_n \geq t)  \leq p^2e^{-nI_1(t/2)} + 2p^2e^{-nI_2(1/2)} \leq 3p^2 e^{-n g(t)}
\eeaa
where $g(t) =\min\{I_1(t/2),\, I_2(1/2)\}$ for any $t>0.$ Obviously,  $g(t)>0$ for any $t>0$ from (i) and (ii) above.\ \ \ \ \ \ $\blacksquare$
\begin{lemma}\lbl{solar.system} Let $Z$ be a random variable with $EZ=0,\, EZ^2=1$ and $Ee^{t_0|Z|}< \infty$ for some $t_0>0.$ Choose $\alpha>0$ such that $E(Z^2e^{\alpha |Z|}) \leq 3/2.$ Set $I(x)=\sup_{t\in \mathbb{R}}\{tx-\log E e^{t Z}\}.$ Then $I(x)\geq x^2/3$ for all $0\leq x \leq 3\alpha/2.$
\end{lemma}
\textbf{Proof}. By the Taylor expansion, for any $x \in \mathbb{R}$, $e^x=1 +x+\frac{x^2}{2}e^
{\theta x}$ for some  $\theta\in [0, 1].$ It follows from $EZ=0$ that
\beaa
Ee^{tZ}=1+\frac{t^2}{2} E(Z^2e^{\theta t Z})\leq 1+ \frac{t^2}{2} E(Z^2e^{t|Z|}) \leq 1+ \frac{3}{4}t^2
\eeaa
for all $0\leq t \leq \alpha.$ Use the inequality $\log (1+x) \leq x$ for all $x>-1$ to see that $\log Ee^{tZ}\leq 3t^2/4$ for every $0\leq t \leq \alpha.$ Take $t_0=2x/3$ with $x>0.$ Then $0\leq t_0 \leq \alpha$ for all $0\leq x \leq 3\alpha/2.$ It follows that
\beaa
I(x) \geq t_0x -\frac{3}{4}t_0^2=\frac{x^2}{3}. \ \ \ \ \ \ \ \ \ \ \ \ \ \ \ \ \ \ \ \ \ \blacksquare
\eeaa

\subsection{Proof of Theorem \ref{birthday}}\label{Limiting}

\begin{lemma}\lbl{long} Let $\xi_1, \cdots, \xi_n$ be i.i.d. random variables with $E\xi_1=0,\, E\xi_1^2=1$ and $Ee^{t_0|\xi_1|^{\alpha}}< \infty$ for some $t_0>0$ and $0<\alpha\leq 1.$ Put $S_n=\sum_{i=1}^n\xi_i$ and $\beta=\alpha/(2+\alpha).$ Then, for any  $\{p_n;\, n\geq 1\}$ with $0<p_n\to \infty$ and $\log p_n =o(n^{\beta})$ and $\{y_n;\, n\geq 1\}$ with $y_n\to y>0,$
\begin{eqnarray*}
P\Big(\frac{S_n}{\sqrt{n\log p_n}}\geq y_n\Big) \sim \frac{p_n^{-y_n^2/2}(\log p_n)^{-1/2}}{\sqrt{2\pi}\, y}
\end{eqnarray*}
as $n\to\infty.$
\end{lemma}

\begin{prop}\lbl{leave} Let $\{x_{ij};\, i\geq 1, j\geq 1\}$ be  i.i.d. random variables with  $Ex_{11}=0$,  $E(x_{11}^2)=1$ and $Ee^{t_0|x_{11}|^{\alpha}}< \infty$ for some $0<\alpha\leq 2$ and $t_0>0.$ Set $\beta=\alpha/(4+\alpha).$  Assume  $p=p(n)\to \infty$ and $\log p=o(n^{\beta})$ as $n\to\infty.$  Then
\beaa
P\left(\frac{W_n^2-\alpha_n}{n} \leq z\right) \to e^{-Ke^{-z/2}}
\eeaa
as $n \to \infty$ for any $z \in \mathbb{R},$ where $\alpha_n=4n\log p-n\log(\log p)$ and $K=(\sqrt{8\pi})^{-1}.$
\end{prop}
\textbf{Proof.} It suffices to show that
\bea\lbl{hartford}
P\Big(\max_{1\leq i < j \leq p}|y_{ij}| \leq \sqrt{\alpha_n+nz}\Big) \to e^{-Ke^{-z/2}},
\eea
where $y_{ij}=\sum_{k=1}^nx_{ki}x_{kj}.$ We now apply Lemma \ref{stein} to prove (\ref{hartford}). Take $I=\{(i,j); 1\leq i < j \leq p\}.$ For $u=(i,j)\in I,$ set $X_{u}=|y_{ij}|$ and $B_{u}=\{(k,l)\in I;\ \mbox{one of}\ k\ \mbox{and}\ l =i\ \mbox{or}\ j,\ \mbox{but}\ (k,l)\ne u\}.$ Let $a_n=\sqrt{\alpha_n+nz}$ and $A_{ij}=\{|y_{ij}|>a_n\}.$ Since $\{y_{ij};\, (i,j)\in I\}$ are identically distributed, by Lemma \ref{stein},
\bea\lbl{protest}
|P(W_n\leq a_n) - e^{-\lambda_n}|\leq  b_{1,n}+b_{2,n}
\eea
where
\bea\lbl{fee}
\lambda_n=\frac{p(p-1)}{2}P(A_{12}),\ b_{1,n}\leq 2p^3P(A_{12})^2\ \mbox{and}\ b_{2,n} \leq 2p^3P(A_{12}A_{13}).
\eea
We first calculate $\lambda_n.$ Write
\bea\lbl{mother}
\lambda_n=\frac{p^2-p}{2}P\left(\frac{|y_{12}|}{\sqrt{n}} > \sqrt{\frac{\alpha_n}{n}+z}\, \right)
\eea
 and $y_{12}=\sum_{i=1}^n\xi_i,$ where $\{\xi_i;\, 1\leq i \leq n\}$ are i.i.d. random variables with the same distribution as  that of $x_{11}x_{12}.$ In particular, $E\xi_1=0$ and $E\xi_1^2=1.$ Note $\alpha_1:=\alpha/2\leq 1.$  We then have
\begin{eqnarray*}
|x_{11}x_{12}|^{\alpha_1}
\leq \left(\frac{x_{11}^2 + x_{12}^2}{2}\right)^{\alpha_1} \leq \frac{1}{2^{\alpha_1}}\Big(|x_{11}|^{\alpha} + |x_{12}|^{\alpha}\Big).
\end{eqnarray*}
Hence, by independence,
\begin{eqnarray*}
 Ee^{t_0|\xi_1|^{\alpha_1}}=Ee^{t_0|x_{11}x_{12}|^{\alpha_1}}<\infty.
\end{eqnarray*}
Let $y_n=\sqrt{(\frac{\alpha_n}{n}+z)/\log p}.$ Then $y_n\to 2$ as $n\to\infty.$ By Lemma \ref{long},
\begin{eqnarray*}
P\left(\frac{y_{12}}{\sqrt{n}} > \sqrt{\frac{\alpha_n}{n}+z}\, \right)
& = & P\left(\frac{\sum_{i=1}^n\xi_i}{\sqrt{n\log p}} > y_n \right)\\
& \sim & \frac{p^{-y_n^2/2}(\log p)^{-1/2}}{2\sqrt{2\pi}}\sim \frac{e^{-z/2}}{\sqrt{8\pi}}\cdot\frac{1}{p^2}
\end{eqnarray*}
as $n\to\infty.$ Considering $Ex_{ij}=0,$ it is easy to see that the above also holds if $y_{12}$ is replaced by $-y_{12}.$ These and (\ref{mother}) imply that
\begin{eqnarray}\lbl{has}
\lambda_n \sim \frac{p^2-p}{2}\cdot 2\cdot \frac{e^{-z/2}}{\sqrt{8\pi}}\cdot\frac{1}{p^2}\sim \frac{e^{-z/2}}{\sqrt{8\pi}}
\end{eqnarray}
as $n\to\infty.$

Recall (\ref{protest}) and (\ref{fee}), to complete the proof,  we have to  verify that $b_{1,n} \to 0$ and $b_{2,n} \to 0$ as $n\to \infty.$  By (\ref{fee}),  (\ref{mother}) and (\ref{has}),
\beaa
b_{1,n}& \leq & 2p^3P(A_{12})^2\\
& = & \frac{8p^3\lambda_n^2}{(p^2-p)^2}=O\left(\frac{1}{p}\right)
\eeaa
as $n\to\infty.$ Also, by (\ref{fee}),
\begin{eqnarray*}
b_{2,n} &\leq & 2p^3P\Big(|y_{12}| > \sqrt{\alpha_n+nz},\ |y_{13}| > \sqrt{\alpha_n+nz}\,\Big)\\
& =& 2p^3 E\Big\{ P^1\Big(|\sum_{k=1}^nx_{k1}x_{k2}| > t_n\sqrt{n\log p}\,\Big)^2\Big\}.
\end{eqnarray*}
where $P^1$ stands for the conditional probability given $\{x_{k,1};\, 1\leq k \leq n\}$, and $t_n: =\sqrt{\alpha_n+nz}/\sqrt{n\log p} \to 2.$ By Lemma \ref{sweet1}, the above expectation is equal to  $O(p^{\epsilon-4})$ as $n\to\infty$ for any $\epsilon>0.$ Now choose $\epsilon \in (0, 1),$ then $b_{2,n}=O(p^{\epsilon-1}) \to 0$ as $n\to\infty.$ The proof is then completed.  \ \ \ \ \ \ $\blacksquare$\\

\noindent\textbf{Proof of Theorem \ref{birthday}}. By the first paragraph in the proof of Theorem \ref{army}, w.l.o.g., assume $\mu=0$ and $\sigma=1.$ From Proposition \ref{leave} and the Slusky lemma, it suffices to show
\begin{eqnarray}\lbl{preschool}
\frac{n^2L_n^2-W_n^2}{n} \to 0
\end{eqnarray}
in probability as $n\to\infty.$ Let $\Delta_n=|nL_n-W_n|$ for $n\geq 1.$ Observe that
\begin{eqnarray}\lbl{cosco}
|n^2L_n^2-W_n^2|=|nL_n-W_n|\cdot |nL_n + W_n| \leq \Delta_n\cdot (\Delta_n + 2W_n).
\end{eqnarray}
It is easy to see from Proposition \ref{leave} that
\begin{eqnarray}\lbl{two}
\frac{W_n}{\sqrt{n\log p}} \to 2
\end{eqnarray}
in probability as $n\to\infty.$ By Lemma \ref{chicago},
\begin{eqnarray*}
\Delta_n \leq \mn n\Gamma_n-X_n^TX_n\mn \leq (b_{n,1}^2+2b_{n,1})W_nb_{n,3}^{-2}+nb_{n,3}^{-2}b_{n,4}^2.
\end{eqnarray*}
By (ii) of Lemma \ref{palo}, $b_{n,3} \to 1$ in probability as $n\to\infty,\ \{\sqrt{n/\log p}\, b_{n,1}\}$ and $\{\sqrt{n/\log p}\, b_{n,4}\}$ are  tight. Set $b_{n,1}'=\sqrt{n/\log p}\, b_{n,1}$ and $b_{n,4}'=\sqrt{n/\log p}\, b_{n,4}$ for all $n\geq 1.$ Then $\{b_{n,1}'\}$ and $\{b_{n,4}'\}$ are  tight. It follows that
\begin{eqnarray*}
\frac{\Delta_n}{\log p}\leq
\left(\sqrt{\frac{\log p}{n}}b_{n,1}'^2 + 2b_{n,1}'\right)\cdot\frac{W_n}{\sqrt{n\log p}}\cdot b_{n,3}^{-2} + b_{n,3}^{-2}b_{n,4}'^2
\end{eqnarray*}
which combining with (\ref{two}) yields that
\begin{eqnarray}\lbl{deed}
\Big\{\frac{\Delta_n}{\log p}\Big\}\ \mbox{is tight}.
\end{eqnarray}
This and (\ref{two}) imply that $\{\Delta_n'\}$ and $\{W_n'\}$ are tight, where $\Delta_n':=\Delta_n/\log p$ and $W_n':=W_n/\sqrt{n\log p}.$ From (\ref{cosco}) and then (\ref{tight}),
\begin{eqnarray}
\frac{|n^2L_n^2-W_n^2|}{n}
&\leq & \frac{(\log p)\Delta_n'\Big\{(\log p)\Delta_n' + 2\sqrt{n\log p}\, W_n'\Big\}}{n}\nonumber\\
& \leq & 2\sqrt{\frac{(\log p)^3}{n}}\Big(\sqrt{\frac{\log p}{n}}\Delta_n' + W_n'\Big)\to 0 \lbl{pay}
\end{eqnarray}
in probability as $n\to\infty$ since $\log p=o(n^{1/3}).$ This gives (\ref{preschool}).\ \ \ \ \ \ \ $\blacksquare$

\subsection{Proof of Theorem \ref{youngth}}\lbl{substantiald}

We begin to prove the Theorem \ref{youngth} by stating three technical
lemmas which are proved in the Appendix.
\begin{lemma}\lbl{webster} Let $\{(u_{k1}, u_{k2}, u_{k3}, u_{k4})^T;\, 1\leq i \leq n\}$ be a sequence of i.i.d. random vectors with distribution $N_4(0, \Sigma_4)$ where
\begin{eqnarray*}
\Sigma_4=\begin{pmatrix}
1& 0 & r & 0\\
0 & 1 & 0 & 0\\
r & 0 & 1 & 0\\
0 & 0 & 0 & 1
\end{pmatrix}
,\ |r|\leq 1.
\end{eqnarray*}
Set $a_n=(4n\log p-n\log(\log p)+ny)^{1/2}$ for $n\geq e^e$ and $y\in \mathbb{R}.$ Suppose $n\to\infty,\, p\to\infty$ with $\log p=o(n^{1/3}).$ Then,
\begin{eqnarray}\lbl{itis}
\sup_{|r|\leq 1}P\Big(|\sum_{k=1}^n u_{k1}u_{k2}|> a_n,\, |\sum_{k=1}^n u_{k3}u_{k4}|> a_n\Big)=O\Big(\frac{1}{p^{4-\epsilon}}\Big)
\end{eqnarray}
for any $\epsilon>0.$
\end{lemma}

\begin{lemma}\lbl{yahoo} Let $\{(u_{k1}, u_{k2}, u_{k3}, u_{k4})^T;\, 1\leq i \leq n\}$ be a sequence of i.i.d. random vectors with distribution $N_4(0, \Sigma_4)$ where
\begin{eqnarray*}
\Sigma_4=\begin{pmatrix}
1& 0 & r_1 & 0\\
0 & 1 & r_2 & 0\\
r_1 & r_2 & 1 & 0\\
0 & 0 & 0 & 1
\end{pmatrix}
,\ \ |r_1|\leq 1,\, |r_2|\leq 1.
\end{eqnarray*}
Set $a_n=(4n\log p-n\log(\log p)+ny)^{1/2}$ for $n\geq e^e$ and $y\in \mathbb{R}.$ Suppose $n\to\infty,\, p\to\infty$ with $\log p=o(n^{1/3}).$ Then, as $n\to\infty,$
\begin{eqnarray*}
\sup_{|r_1|,\, |r_2|\leq 1}P\Big(|\sum_{k=1}^n u_{k1}u_{k2}|> a_n,\, |\sum_{k=1}^n u_{k3}u_{k4}|> a_n\Big)=O\Big(p^{-\frac{8}{3}+\epsilon}\Big)
\end{eqnarray*}
for any $\epsilon>0.$
\end{lemma}

\begin{lemma}\lbl{yahoo1} Let $\{(u_{k1}, u_{k2}, u_{k3}, u_{k4})^T;\, 1\leq i \leq n\}$ be a sequence of i.i.d. random vectors with distribution $N_4(0, \Sigma_4)$ where
\begin{eqnarray*}
\Sigma_4=\begin{pmatrix}
1& 0 & r_1 & 0\\
0 & 1 & 0 & r_2\\
r_1 & 0 & 1 & 0\\
0 & r_2 & 0 & 1
\end{pmatrix}
,\ \ |r_1|\leq 1,\, |r_2|\leq 1.
\end{eqnarray*}
Set $a_n=(4n\log p-n\log(\log p)+ny)^{1/2}$ for $n\geq e^e$ and $y\in \mathbb{R}.$ Suppose $n\to\infty,\, p\to\infty$ with $\log p=o(n^{1/3}).$ Then, for any $\delta\in (0, 1),$ there exists $\epsilon_0=\epsilon(\delta)>0$ such that
\begin{eqnarray}\lbl{demand}
\sup_{|r_1|,\, |r_2|\leq 1-\delta}P\Big(|\sum_{k=1}^n u_{k1}u_{k2}|> a_n,\, |\sum_{k=1}^n u_{k3}u_{k4}|> a_n\Big)=O\Big(p^{-2-\epsilon_0}\Big)
\end{eqnarray}
as $n\to\infty.$
\end{lemma}

Recall notation $\tau$, $\Sigma=(\sigma_{ij})_{p\times p}$ and $X_n=(x_{ij})_{n\times p}\sim N_p(\mu, \Sigma)$ above (\ref{Hypo}).

\begin{prop}\lbl{correlated} Assume $\mu=0$ and $\sigma_{ii}=1$ for all $1\leq i \leq p.$ Define
\begin{eqnarray}\lbl{jinmen}
V_n=V_{n,\tau}=\max_{1\leq i < j\leq p,\, |j-i|\geq \tau}|x_i^Tx_j|.
\end{eqnarray}
Suppose $n\to\infty,\, p=p_n\to\infty$ with $\log p=o(n^{1/3})$, $\tau=o(p^t)$ for any $t>0$, and for some $\delta\in (0,1)$, $|\Gamma_{p,\delta}|=o(p)$ as $n\to\infty$.  Then, under $H_0$ in (\ref{Hypo}),
\beaa
P\left(\frac{V_n^2-\alpha_n}{n} \leq y\right) \to e^{-Ke^{-y/2}}
\eeaa
as $n \to \infty$ for any $y \in \mathbb{R},$ where $\alpha_n=4n\log p-n\log(\log p)$ and $K=(\sqrt{8\pi})^{-1}.$
\end{prop}
\noindent\textbf{Proof}. Set $a_n=(4n\log p-n\log(\log p) + ny)^{1/2},$
\begin{eqnarray}
& &  \Lambda_{p}=\Big\{(i, j);\ 1\leq i < j \leq p,\ j-i\geq \tau,\, \max_{1\leq k\ne i \leq p}\{|r_{ik}|\} \leq  1-\delta,\, \max_{1\leq k\ne j \leq p}\{|r_{jk}|\} \leq  1-\delta\Big\}, \nonumber\\
& & V_n'=\max_{(i,j)\in \Lambda_{p}}\Big|\sum_{k=1}^nx_{ki}x_{kj}\Big|.\lbl{aobao}
\end{eqnarray}
\noindent{\it Step 1}. We claim that, to prove the proposition, it suffices to show
\begin{eqnarray}\lbl{August}
\lim _{n\to\infty}P\left(V_n'\leq  a_n\right) = e^{-Ke^{-y/2}}
\end{eqnarray}
for any $y\in\mathbb{R}.$

In fact, to prove the theorem, we need to show that
\begin{eqnarray}\lbl{joking}
\lim _{n\to\infty}P\left(V_n> a_n\right) = 1-e^{-Ke^{-y/2}}
\end{eqnarray}
for every $y\in \mathbb{R}.$ Notice $\{x_{ki},\, x_{kj};\, 1\leq k \leq n\}$ are $2n$ i.i.d. standard normals if $|j-i|\geq \tau.$ Then
\begin{eqnarray*}
P\left(V_n> a_n\right)&\leq & P\left(V_n'> a_n\right) + \sum P\Big(|\sum_{k=1}^nx_{k1}x_{k\, \tau+1}|> a_n\Big)
\end{eqnarray*}
where the sum runs over all pair $(i,j)$  such that $1\leq i<j \leq p$ and one of $i$ and $j$ is in $\Gamma_{p, \delta}.$
Note that $|x_{11}x_{1\, \tau+1}|\leq (x_{11}^2 + x_{1\, \tau+1}^2)/2,$ it follows that $Ee^{|x_{11}x_{1\, \tau+1}|/2}<\infty$ by independence and $E\exp(N(0,1)^2/4)<\infty.$ Since $\{x_{k1}, x_{k\, \tau+1};\, 1\leq k \leq n\}$ are i.i.d. with mean zero and variance one, and $y_n:=a_n/\sqrt{n\log p}\to 2$ as $n\to\infty,$ taking $\alpha=1$ in Lemma \ref{long}, we get
\begin{eqnarray}
& & P\Big(\frac{1}{\sqrt{n\log p}}|\sum_{k=1}^nx_{k1}x_{k\, \tau+1}|> \frac{a_n}{\sqrt{n\log p}}\Big)\nonumber\\
& \sim & 2\cdot \frac{p^{-y_n^2/2}(\log p)^{-1/2}}{2\sqrt{2\pi}}\sim \frac{e^{-y/2}}{\sqrt{2\pi}}\cdot\frac{1}{p^2}\lbl{Gooden}
\end{eqnarray}
as $n\to\infty.$ Moreover, note that the total number of such pairs is no more than $2p\,|\Gamma_{p,\delta}|.$ Therefore,
\begin{eqnarray}
P\left(V_n'> a_n\right) \leq  P\left(V_n> a_n\right)
& \leq & P\left(V_n'> a_n\right) + 2p\,|\Gamma_{p,\delta}|\cdot P\Big(|\sum_{k=1}^nx_{k1}x_{k\, \tau+1}|> a_n\Big)\nonumber\\
& \leq & P\left(V_n'> a_n\right) + o(p^2)\cdot O\Big(\frac{1}{p^2}\Big)\lbl{gongbing}
\end{eqnarray}
by the assumption on $\Gamma_{p,\delta}$ and (\ref{Gooden}). Thus, this joint with (\ref{joking}) gives (\ref{August}).

\noindent{\it Step 2}. We now apply Lemma \ref{stein} to prove (\ref{August}). Take $I=\Lambda_p.$ For $(i,j)\in I,$ set $Z_{ij}=|\sum_{k=1}^nx_{ki}x_{kj}|$,
\begin{eqnarray*}
& &  B_{i,j}=\{(k,l)\in \Lambda_p;\ |s-t| < \tau\ \mbox{for some}\ s\in \{k, l\}\ \mbox{and some}\ t\in\{i, j\},\,   \mbox{but}\ (k,l)\ne (i,j)\},\ \ \ \ \ \ \ \ \ \ \ \\
& & a_n=\sqrt{\alpha_n+ny}\ \ \mbox{and}\ \  A_{ij}=\{|Z_{ij}|>a_n\}.
\end{eqnarray*}
It is easy to see that $|B_{i,j}|\leq 2\cdot (2\tau + 2\tau)p=8\tau p$ and that $Z_{ij}$ are independent of $\{Z_{kl};\,(k,l) \in \Lambda_p\backslash B_{i,j}\}$ for any $(i,j)\in \Lambda_p.$ By Lemma \ref{stein},
\bea\lbl{protest1}
|P(V_n\leq a_n) - e^{-\lambda_n}|\leq  b_{1,n}+b_{2,n}
\eea
where
\bea
& & \lambda_n=|\Lambda_p|\cdot P(A_{1\, \tau+1}),\ b_{1,n}\leq \sum_{d \in \Lambda_p}\sum_{d' \in B_{a}}P(A_{12})^2=8\tau p^3P(A_{1\, \tau+1})^2\ \ \mbox{and}\lbl{elk1}\ \ \ \ \ \ \ \ \ \ \ \  \ \\\
& & b_{2,n} \leq \sum_{d \in \Lambda_p}\sum_{d\ne d' \in B_{a}}P(Z_{d} >t, Z_{d'} >t)\lbl{elk2}\ \ \ \ \ \ \ \ \ \ \ \  \ \ \ \ \ \ \ \ \ \ \  \ \ \ \ \ \ \ \ \ \ \
\eea
 from the fact that $\{Z_{ij};\, (i,j)\in \Lambda_p\}$ are identically distributed. We first calculate $\lambda_n.$  By definition
\begin{eqnarray*}
\frac{p^2}{2} > |\Lambda_p|& \geq & \Big|\{(i, j);\ 1\leq i < j \leq p,\ j-i\geq \tau\}\Big| - 2p\cdot |\Gamma_{p,\delta}|\\
& = & \sum_{i=1}^{p-\tau}(p-\tau-i+1) - 2p\cdot |\Gamma_{p,\delta}|.
\end{eqnarray*}
Now the sum above is equal to $\sum_{j=1}^{p-\tau}j=(p-\tau)(p-\tau+1)/2 \sim p^2/2$ since $\tau=o(p)$. By assumption  $|\Gamma_{p,\delta}|=o(p)$ we conclude that
\begin{eqnarray}\lbl{flooring}
|\Lambda_p| \sim \frac{p^2}{2}
\end{eqnarray}
as $n\to\infty.$ It then follows from  (\ref{Gooden}) that
\bea\lbl{mother1}
\lambda_n \sim \frac{p^2}{2}\cdot \frac{e^{-y/2}}{\sqrt{2\pi}}\cdot\frac{1}{p^2}\sim \frac{e^{-y/2}}{\sqrt{8\pi}}
\eea
as $n\to\infty.$

Recall (\ref{protest1}) and (\ref{mother1}), to complete the proof,  we have to  verify that $b_{1,n} \to 0$ and $b_{2,n} \to 0$ as $n\to \infty.$ Clearly,  by the first expression in (\ref{elk1}), we get from (\ref{mother1}) and then (\ref{flooring}) that
\beaa
b_{1,n}& \leq & 8\tau p^3P(A_{1\, \tau+1})^2 =  \frac{8\tau p^3\lambda_n^2}{|\Lambda_{p}|^2}=O\left(\frac{\tau}{p}\right)\to 0
\eeaa
as $n\to\infty$ by the assumption on $\tau.$

\noindent{\it Step 3}. Now we consider $b_{2,n}.$ Write $d=(d_1, d_2)\in \Lambda_p$ and $d'=(d_3, d_4)\in \Lambda_p$ with $d_1<d_2$ and $d_3< d_4.$ It is easy to see from (\ref{elk2}) that
\begin{eqnarray*}
b_{2,n}\leq 2\sum P(Z_{d} >a_n, Z_{d'} > a_n)
\end{eqnarray*}
where the sum runs over every pair $(d, d')$ satisfying
\begin{eqnarray}\lbl{helicopter}
d, d'\in \Lambda_p,\ d\ne d',\ d_1\leq d_3\ \mbox{and}\ |d_i-d_j|< \tau\ \mbox{for some}\ i\in \{1,2\}\ \mbox{and some}\ j\in\{3,4\}.
\end{eqnarray}
Geometrically, there are three cases for the locations of $d=(d_1, d_2)$ and $d'=(d_3, d_4)$:
\begin{eqnarray}\lbl{fire}
(1)\, d_2 \leq d_3;\ \ (2)\, d_1\leq d_3< d_4 \leq d_2;\ \ (3)\, d_1\leq d_3\leq d_2 \leq d_4.
\end{eqnarray}
Let $\Omega_j$ be the subset of index $(d, d')$ with restrictions (\ref{helicopter}) and $(j)$ for $j=1,2,3.$  Then
\begin{eqnarray}\lbl{frog}
b_{2,n}\leq 2\sum_{i=1}^3\sum_{(d, d')\in \Omega_i} P(Z_{d} >a_n, Z_{d'} > a_n).
\end{eqnarray}
We next analyze each of the three sums separately. Recall all diagonal entries of $\Sigma$ in $N_p(0,\, \Sigma)$ are equal to $1.$ Let random vector
\begin{eqnarray}\lbl{coordinate}
(w_1, w_2, \cdots, w_p) \sim N_p(0,\, \Sigma).\ \ \ \ \ \ \ \ \ \ \ \ \ \ \ \ \ \ \ \ \ \ \ \ \ \ \ \
\end{eqnarray}
Then every $w_i$ has the distribution of $N(0, 1).$

\noindent\textbf{Case (1)}. Evidently, (\ref{helicopter}) and (1) of (\ref{fire}) imply that $0\leq d_3- d_2 < \tau.$ Hence,  $|\Omega_1| \leq \tau p^3.$ Further, for $(d, d')\in \Omega_1,$ the covariance matrix of $(w_{d_1}, w_{d_2}, w_{d_3}, w_{d_4})$ is equal to
\begin{eqnarray*}
\begin{pmatrix}
1& 0 & 0 & 0\\
0 & 1 & \gamma & 0\\
0 & \gamma & 1 & 0\\
0 & 0 & 0 & 1
\end{pmatrix}
\end{eqnarray*}
for some $\gamma\in [-1,1].$ Thus, the covariance matrix of $(w_{d_2}, w_{d_1}, w_{d_3}, w_{d_4})$ is equal to
\begin{eqnarray*}
\begin{pmatrix}
1& 0 & \gamma & 0\\
0 & 1 & 0 & 0\\
\gamma & 0 & 1 & 0\\
0 & 0 & 0 & 1
\end{pmatrix}
.
\end{eqnarray*}
Recall $Z_d=Z_{d_1,d_2}=Z_{d_2, d_1}=|\sum_{k=1}^nx_{kd_1}x_{kd_2}|$ defined at the beginning of {\it Step 2}. By Lemma \ref{webster}, for some $\epsilon>0$ small enough,
\begin{eqnarray}\lbl{Fourier}
\sum_{(d, d')\in \Omega_1} P(Z_{d} >a_n, Z_{d'} > a_n)\nonumber
&= &\sum_{(d, d')\in \Omega_1} P(Z_{d_2, d_1} >a_n, Z_{d_3, d_4} > a_n)\nonumber\\
& \leq & \tau p^3\cdot O\Big(\frac{1}{p^{4-\epsilon}}\Big)= O\Big(\frac{\tau}{p^{1-\epsilon}}\Big)\to 0
\end{eqnarray}
as $n\to\infty$ since $\tau=o(p^{t})$ for any $t>0.$

\noindent\textbf{Case (2)}. For any $(d, d')\in \Omega_2,$ there are three possibilities.

\noindent (I): $|d_1-d_3|< \tau$ and $|d_2-d_4|< \tau$;\, (II): $|d_1-d_3|< \tau$ and $|d_2-d_4|\geq \tau$;\, (III): $|d_1-d_3|\geq  \tau$ and $|d_2-d_4|< \tau$. The case that $|d_1-d_3|\geq  \tau$ and $|d_2-d_4| \geq \tau$ is excluded by   (\ref{helicopter}).

 Let $\Omega_{2, I}$ be the subset of $(d, d')\in \Omega_2$ satisfying (I), and $\Omega_{2, II}$ and $\Omega_{2, III}$ be  defined similarly. It is easy to check that $|\Omega_{2, I}| \leq \tau^2p^2.$ The covariance matrix of $(w_{d_1}, w_{d_2}, w_{d_3}, w_{d_4})$  is equal to
\begin{eqnarray*}
\begin{pmatrix}
1& 0 & \gamma_1 & 0\\
0 & 1 & 0 & \gamma_2\\
\gamma_1 & 0 & 1 & 0\\
0 & \gamma_2 & 0 & 1
\end{pmatrix}
\end{eqnarray*}
for some $\gamma_1, \gamma_2\in [-1,1].$ By Lemma \ref{yahoo1},
\begin{eqnarray}\lbl{I}
\sum_{(d, d')\in \Omega_{2, I}} P(Z_{d} > a_n, Z_{d'} > a_n)= O\Big(\frac{\tau^2}{p^{\epsilon_0}}\Big)\to 0
\end{eqnarray}
as $n\to\infty.$

Observe $|\Omega_{2, II}| \leq \tau p^3.$ The covariance matrix of $(w_{d_1}, w_{d_2}, w_{d_3}, w_{d_4})$  is equal to
\begin{eqnarray*}
\begin{pmatrix}
1& 0 & \gamma & 0\\
0 & 1 & 0 & 0\\
\gamma & 0 & 1 & 0\\
0 & 0 & 0 & 1
\end{pmatrix}
,\ \ |\gamma|\leq 1.
\end{eqnarray*}
By Lemma \ref{webster}, take $\epsilon>0$ small enough to get
\begin{eqnarray}\lbl{II}
\sum_{(d, d')\in \Omega_{2, II}} P(Z_{d} > a_n, Z_{d'} > a_n)= O\Big(\frac{\tau}{p^{1-\epsilon}}\Big)\to 0
\end{eqnarray}
as $n\to\infty.$

The third case is similar to the second one. In fact,  $|\Omega_{2, III}| \leq \tau p^3.$ The covariance matrix of $(w_{d_1}, w_{d_2}, w_{d_3}, w_{d_4})$  is equal to
\begin{eqnarray*}
\begin{pmatrix}
1& 0 & 0 & 0\\
0 & 1 & 0 & \gamma\\
0 & 0 & 1 & 0\\
0 & \gamma & 0 & 1
\end{pmatrix}
,\ \ |\gamma|\leq 1.
\end{eqnarray*}
Thus, the covariance matrix of $(w_{d_2}, w_{d_1}, w_{d_4}, w_{d_3})$  is equal to $\Sigma_4$ in Lemma \ref{webster}. Then, by the same argument as that in the equality in (\ref{Fourier}) we get
\begin{eqnarray}\lbl{III}
\sum_{(d, d')\in \Omega_{2, III}} P(Z_{d} > a_n, Z_{d'} > a_n)= O\Big(\frac{\tau}{p^{1-\epsilon}}\Big)\to 0
\end{eqnarray}
as $n\to\infty$ by taking $\epsilon>0$ small enough. Combining (\ref{I}), (\ref{II}) and (\ref{III}), we conclude
\begin{eqnarray*}
\sum_{(d, d')\in \Omega_2} P(Z_{d} >a_n, Z_{d'} > a_n)\to 0
\end{eqnarray*}
as $n\to\infty.$ This and (\ref{Fourier}) together with (\ref{frog}) say that, to finish the proof of this proposition, it suffices to verify
\begin{eqnarray}\lbl{Suzuki}
\sum_{(d, d')\in \Omega_3} P(Z_{d} >a_n, Z_{d'} > a_n)\to 0
\end{eqnarray}
as $n\to\infty.$ The next lemma confirms this. The proof is then completed.\ \ \ \ \ \ \ \ $\blacksquare$

\begin{lemma}\lbl{fertilizer} Let the notation be as in the proof of Proposition \ref{correlated}, then (\ref{Suzuki}) holds.
\end{lemma}

\noindent\textbf{Proof of Theorem \ref{youngth}}. By the first paragraph in the proof of Theorem \ref{army}, w.l.o.g., we prove the theorem  by assuming that the $n$ rows of $X_n=(x_{ij})_{1\leq i \leq n, 1\leq j \leq p}$ are i.i.d. random vectors with distribution $N_p(0, \Sigma)$ where all of the diagonal entries of $\Sigma$ are equal to $1.$ Consequently,  by the assumption on $\Sigma,$ for any subset $E=\{i_1, i_2, \cdots, i_m\}$ of $\{1, 2, \cdots, p\}$ with $|i_s-i_t| \geq \tau$ for all $1\leq s< t \leq m,$ we know that $\{x_{ki};\, 1\leq k \leq n,\, i\in E\}$ are $mn$ i.i.d. $N(0, 1)$-distributed random variables.

Reviewing the proof of Lemma \ref{palo}, the argument is only based on the distribution of each column of $\{x_{ij}\}_{n\times p};$ the joint distribution of any two different columns are irrelevant. In current situation, the entries in each column are i.i.d. standard normals. Thus, take $\alpha=2$ in the lemma to have
\begin{eqnarray}\lbl{cloud}
& & b_{n,3}\to 1\ \mbox{in probability as}\ n\to\infty, \nonumber\\
& & \Big\{\sqrt{\frac{n}{\log p}}\, b_{n,1}\Big\}\ \mbox{and}\ \Big\{\sqrt{\frac{n}{\log p}}\, b_{n,4}\Big\}\ \mbox{are  tight}\ \ \ \ \ \ \ \ \ \ \ \ \ \ \ \ \ \ \ \ \ \ \ \ \ \ \ \
\end{eqnarray}
as $n\to\infty$, $p\to \infty$ with $\log p=o(n),$ where $b_{n,1},\, b_{n,3}$ and $b_{n,4}$ are as in Lemma \ref{palo}. Let $V_n=V_{n, \tau}=(v_{ij})_{p\times p}$ be as in (\ref{jinmen}). It is seen from Proposition \ref{correlated} that
\bea\lbl{sunflower}
\frac{V_{n,\tau}}{\sqrt{n\log p}} \to 2
\eea
in probability as $n\to\infty$, $p\to\infty$ and $\log p=o(n^{1/3}).$ Noticing the differences in the indices of $\max_{1\leq i <  j \leq p}|\rho_{ij}|$ and $\max_{1\leq i <  j \leq p,\, |i-j|\geq \tau}|\rho_{ij}|=L_{n,\tau}$, checking the proof of Lemma 2.2 from Jiang (2004a), it is easy to see that
\bea\lbl{county}
\Delta_n:=\max_{1\leq i <  j \leq p,\, |i-j|\geq \tau}\Big|n\rho_{ij}-v_{ij}\Big| \leq (b_{n,1}^2+2b_{n,1})V_{n,\tau}b_{n,3}^{-2}+nb_{n,3}^{-2}b_{n,4}^2.
\eea
Now, using (\ref{cloud}), (\ref{sunflower}) and (\ref{county}), replacing $W_n$ with $V_{n,\tau}$ and $L_n$ with $L_{n,\tau}$ in the proof of Theorem \ref{birthday}, and repeating the whole proof again, we obtain
\begin{eqnarray*}
\frac{n^2L_{n,\tau}^2-V_{n, \tau}^2}{n} \to 0
\end{eqnarray*}
in probability as $n\to\infty.$ This joint with Proposition \ref{correlated} and the Slusky lemma yields the desired limiting result for $L_{n,\tau}.$ \ \ \ \ \ \ \ $\blacksquare$

\newpage
\section{Appendix}\lbl{ladies}

In this appendix we prove Proposition \ref{sister} and verify
the three examples given in Section \ref{CS.sec}. We then prove Lemmas
\ref{palo} - \ref{sweet1} and Lemmas \ref{long} - \ref{fertilizer} which are used in the proof of the main results.

\bigskip
\noindent\textbf{Proof of Proposition \ref{sister}}. We prove the proposition by following the outline of the proof of Proposition \ref{father} step by step. It suffices to show
\begin{eqnarray}
& & \lim_{n\to\infty}P\Big(\frac{W_n}{\sqrt{n\log p}} \geq  2+2\epsilon\Big)= 0\ \ \mbox{and}\ \lbl{winding}\\
& & \lim_{n\to\infty}P\Big(\frac{W_n}{\sqrt{n\log p}} \leq  2-\epsilon\Big)= 0\lbl{crazy}
\end{eqnarray}
for any $\epsilon>0$ small enough. Note that $|x_{11}x_{12}|^{\varrho}=|x_{11}|^{\varrho}\cdot|x_{12}|^{\varrho}\leq |x_{11}|^{2\varrho} +|x_{12}|^{2\varrho}$ for any $\varrho>0.$ The given moment condition implies that $E\exp\big(t_0|x_{11}|^{4\beta/(1-\beta)}\big)< \infty.$ Hence $E\exp\big(|x_{11}|^{\frac{4\beta}{1 + \beta}}\big)< \infty$ and $E\exp\big(|x_{11}x_{12}|^{\frac{2\beta}{1 + \beta}}\Big)< \infty.$ By (i) of Lemma \ref{Xia}, (\ref{sky}) holds for $\{p_n\}$ such that $p_n \to\infty$ and $\log p_n=o(n^{\beta}).$ By using (\ref{grass}) and (\ref{kick}), we obtain (\ref{winding}).

 By using condition $E\exp\{t_0|x_{11}|^{\frac{4\beta}{1 + \beta}}\}< \infty$ again, we know (\ref{hurdle}) also holds for $\{p_n\}$ such that $p_n\to\infty$ and $\log p_n=o(n^{\beta}).$ Then all statements after (\ref{tree}) and before (\ref{sleeping}) hold. Now, by Lemma \ref{sweet1}, (\ref{Nike})  holds for $\{p_n\}$ such that $p_n\to\infty$ and  $\log p_n=o(n^{\beta}),$ we then have (\ref{littleshow}). This implies (\ref{tree}), which is the same as  (\ref{crazy}).\ \ \ \ \ \ \ \ $\blacksquare$\\

\noindent\textbf{Verifications of (\ref{chop}), (\ref{karaok}) and (\ref{roasted})}. We consider the three one by one.

\noindent\textbf{(i)} If $x_{11}\sim N(0, n^{-1})$ as in (\ref{normal.CS}), then $\xi$ and $\eta$ are i.i.d. with distribution $N(0, 1).$ By Lemma 3.2 from Jiang (2005), $I_2(x)=(x-1-\log x)/2$ for $x>0.$ So $I_2(1/2)>1/12.$ Also, since $Ee^{\theta \xi \eta}=Ee^{\theta^2\xi^2/2}=(1-\theta^2)^{-1/2}$ for $|\theta|<1.$ It is straightforward to get
\beaa
I_1(x)=\frac{\sqrt{4x^2+1}-1}{2} - \frac{1}{2}\log \frac{\sqrt{4x^2+1} + 1}{2},\ \  x > 0.
\eeaa
Let $y=\frac{\sqrt{4x^2+1}-1}{2}.$ Then $y>2x^2/3$ for all $|x|\leq 4/5.$ Thus, $I_1(x)=y-\frac{1}{2}\log (1+y)>\frac{y}{2}>\frac{x^2}{3}$ for $|x|\leq 4/5.$ Therefore, $g(t)\geq \min\{I_1(\frac{t}{2}),\, \frac{1}{12}\} \geq \min\{\frac{t^2}{12}, \frac{1}{12}\}=\frac{t^2}{12}$ for $|t| \leq  1.$ Since $1/(2k-1) \leq 1$ if $k\geq 1.$ By Proposition \ref{applesensing}, we have
\bea\lbl{chop.prime}
P\left((2k-1)\tilde{L}_n < 1\right) \geq 1- 3 p^2\exp\Big\{-\frac{n}{12(2k-1)^2}\Big\}
\eea
for all $n\geq 2$ and $k\geq 1,$ which is  (\ref{chop}).

\noindent\textbf{(ii)} Let $x_{11}$ be such that $P(x_{11}=\pm 1/\sqrt{n})=1/2$ as in (\ref{bernoulli.CS}). Then $\xi$ and $\eta$ in Proposition \ref{applesensing} are i.i.d. with $P(\xi=\pm 1)=1/2.$ Hence, $P(\xi\eta=\pm 1)=1/2$ and $\xi^2=1.$ Immediately, $I_2(1)=0$ and $I_2(x)=+\infty$ for all $x \ne 1.$ If $\alpha=\log \frac{3}{2} \sim 0.405,$ then $E(Z^2e^{\alpha |Z|})=e^{\alpha} \leq \frac{3}{2}$ with $Z=\xi\eta.$ Thus, by Lemma \ref{solar.system}, $I_1(x)\geq x^2/3$ for all $0\leq x \leq \frac{3}{5}\leq \frac{3\alpha}{2}.$ Therefore, $g(t) \geq \frac{t^2}{12}$ for $0\leq t \leq \frac{6}{5}.$ This gives that
\bea\lbl{karaok.prime}
P\left((2k-1)\tilde{L}_n < 1\right) \geq 1- 3 p^2\exp\Big\{-\frac{n}{12(2k-1)^2}\Big\}
\eea
provided $\frac{1}{2k-1} \leq \frac{6}{5},$ that is, $k\geq \frac{11}{12}.$ We then obtain (\ref{karaok}) since $k$ is an integer.

\noindent\textbf{(iii)} Let $x_{11}$ be such that $P(x_{11}=\pm \sqrt{3/n})=1/6$ and $P(x_{11}=0)=2/3$ as in (\ref{sparse.CS}). Then $\xi$ and $\eta$ in Proposition \ref{applesensing} are i.i.d. with  $P(\xi=\pm \sqrt{3})=1/6$ and $P(\xi=0)=2/3.$ It follows that $P(Z=\pm 3)=1/18$ and $P(Z=0)=8/9$ with $Z=\xi\eta.$ Take $\alpha=\frac{1}{3}\log \frac{3}{2}> 0.13.$ Then $E(Z^2e^{\alpha |Z|})=\frac{2\times 9}{18}e^{3\alpha}= \frac{3}{2}.$ Thus, by Lemma \ref{solar.system}, $I_1(x)\geq x^2/3$ for all $0\leq x \leq \frac{3\alpha}{2}=\frac{1}{2}\log \frac{3}{2} \sim 0.2027.$ Now, $P(\xi^2=3)=\frac{1}{3}=1-P(\xi^2=0).$ Hence, $\xi^2/3 \sim Ber(p)$ with $p=\frac{1}{3}.$ It follows that
\beaa
I_2(x) &= & \sup_{\theta \in \mathbb{R}}\Big\{(3\theta) \frac{x}{3}-\log Ee^{3\theta (\xi^2/3)}\Big\}\\
& = & \Lambda^*\Big(\frac{x}{3}\Big)=\frac{x}{3}\log x + \Big(1-\frac{x}{3}\Big)\log \frac{3-x}{2}
\eeaa
for $0\leq x \leq 3$ by (b) of Exercise 2.2.23 from \cite{DZ98}. Thus, $I_2(\frac{1}{2})=\frac{1}{6}\log \frac{1}{2} + \frac{5}{6}\log \frac{5}{4}\sim 0.0704>\frac{1}{15}.$ Now, for $0\leq t \leq \frac{2}{5},$ we have
\beaa
g(t)=\min \Big\{I_1(\frac{t}{2}),\, I_2(\frac{1}{2})\Big\} \geq \min\Big\{\frac{t^2}{12},\, \frac{1}{15}\Big\} =\frac{t^2}{12}.
\eeaa
Easily, $t:=\frac{1}{2k-1}\leq \frac{2}{5}$ if and only if $k\geq \frac{7}{4}.$ Thus,  by Proposition \ref{applesensing},
\bea\lbl{roasted.prime}
P\left((2k-1)\tilde{L}_n < 1\right) \geq 1- 3 p^2\exp\Big\{-\frac{n}{12(2k-1)^2}\Big\}
\eea
for all $n\geq 2$ and $k\geq \frac{7}{4}.$ We finally conclude (\ref{roasted}) since $k$ is an integer. \ \ \ \ \ \ \ $\blacksquare$

\bigskip
\noindent\textbf{Proof of Lemma \ref{palo}}. (i) First, since $x_{ij}$'s are i.i.d. bounded random variables with mean zero and variance one, by (i) of Lemma \ref{Xia},
\begin{eqnarray}
P(\sqrt{n/\log p}\, b_{n,4} \geq K)
&= & P\Big(\max_{1\leq i \leq p}\left|\frac{1}{\sqrt{n\log p}}\sum_{k=1}^nx_{ki}\right|\geq K\Big)\lbl{Oregon}\\
& \leq & p\cdot P\Big(\Big|\frac{1}{\sqrt{n\log p}}\sum_{k=1}^nx_{k1}\Big|\geq K\Big)\nonumber\\
& \leq & p\cdot e^{-(K^2/3)\log p}=\frac{1}{p^{K^2/3-1}} \to 0\lbl{Tony}
\end{eqnarray}
as $n\to\infty$ for any $K>\sqrt{3}.$ This says that $\{\sqrt{n/\log p}\, b_{n,4}\}$ are tight.

Second, noticing that $|t-1| \leq |t^2-1|$ for any $t>0$ and  $nh_i^2=\|x_i-\bar{x}_i\|^2=x_i^Tx_i-n|\bar{x}_i|^2,$  we get that
\bea\lbl{baby}
b_{n,1} \leq \max_{1\leq i \leq p}|h_i^2 -1| &\leq & \max_{1\leq i \leq p}\left|\frac{1}{n}\sum_{k=1}^n(x_{ki}^2 -1)\right| + \max_{1\leq i \leq p}\left|\frac{1}{n}\sum_{k=1}^nx_{ki}\right|^2\nonumber\\
& = & Z_n + b_{n,4}^2
\eea
where $Z_n=\max_{1\leq i \leq p}\left|\frac{1}{n}\sum_{k=1}^n(x_{ki}^2 -1)\right|.$
Therefore,
\begin{eqnarray}\lbl{seashore}
\sqrt{\frac{n}{\log p}}b_{n,1} \leq \sqrt{\frac{n}{\log p}}Z_n + \sqrt{\frac{\log p}{n}}\cdot \Big(\sqrt{\frac{n}{\log p}}b_{n,4}\Big)^2.
\end{eqnarray}
Replacing ``$x_{ki}$" in (\ref{Oregon}) with ``$x_{ki}^2-1$" and using the same argument, we obtain that $\{\sqrt{n/\log p}\, Z_{n}\}$ are tight.
Since $\log p =o(n)$ and $\{\sqrt{n/\log p}\, b_{n,4}\}$ are tight, using (\ref{tight}) we know the second term on the right hand side of (\ref{seashore}) goes to zero in probability as $n\to\infty.$ Hence, we conclude from (\ref{seashore}) that $\{\sqrt{n/\log p}\, b_{n,1}\}$ are tight.

Finally, since $\log p =o(n)$ and $\{\sqrt{n/\log p}\, b_{n,1}\}$ are tight, use (\ref{tight}) to have $b_{n,1}\to 0$ in probability as $n\to\infty.$ This implies that $b_{n,3} \to 1$ in probability as $n\to \infty.$

(ii) By (\ref{Tony}) and (\ref{seashore}), to prove the conclusion, it is enough to show, for some constant $K>0,$
\begin{eqnarray}
& & p\cdot P\Big(\,\Big|\frac{1}{\sqrt{n\log p}}\sum_{k=1}^nx_{k1}\Big|\geq K\Big)\to 0\ \ \ \mbox{and}\ \ \ \ \ \ \ \ \ \ \ \ \ \ \ \ \ \ \ \ \ \ \ \ \lbl{sound}\\
& & p\cdot P\Big(\,\Big|\frac{1}{\sqrt{n\log p}}\sum_{k=1}^n(x_{k1}^2-1)\Big|\geq K\Big)\to 0\lbl{swag}
\end{eqnarray}
as $n\to\infty.$ Using $a_n:=\sqrt{\log p_n}=o(n^{\beta/2})$ and (i) of Lemma \ref{Xia},  we have
\begin{eqnarray*}
& & P\Big(\,\Big|\frac{1}{\sqrt{n\log p}}\sum_{k=1}^nx_{k1}\Big|\geq K\Big) \leq \frac{1}{p^{K^2/3}}\ \ \mbox{and}\\
& & P\Big(\,\Big|\frac{1}{\sqrt{n\log p}}\sum_{k=1}^n(x_{k1}^2-1)\Big|\geq K\Big)\leq  \frac{1}{p^{K^2/3}}
\end{eqnarray*}
as $n$ is sufficiently large, where the first inequality holds provided   $E\exp\big(t_0|x_{11}|^{2\beta/(1+\beta)}\big)=E\exp(t_0|x_{11}|^{\alpha/2})<\infty;$ the second holds since  $E\exp\big(t_0|x_{11}^2-1|^{2\beta/(1+\beta)}\big)=E\exp(t_0|x_{11}^2-1|^{\alpha/2})< \infty$ for some $t_0>0,$ which is equivalent to $Ee^{t_0'|x_{11}|^{\alpha}}<\infty$ for some $t_0'>0.$ We then get (\ref{sound}) and (\ref{swag}) by taking $K=2.$ \ \ \ \ \ \ \ $\blacksquare$\\

\noindent\textbf{Proof of Lemma \ref{sweet}}. Let $G_n=\{|\sum_{k=1}^nx_{k1}^2/n-1|< \delta\}.$ Then, by the Chernoff bound (see, e.g., p. 27 from Dembo and Zeitouni (1998)), for any $\delta \in (0,1),$ there exists a constant $C_{\delta}>0$ such that $P(G_n^c)\leq 2e^{-nC_{\delta}}$ for all $n\geq 1.$ Set $a_n=t_n\sqrt{n\log p}.$ Then
\begin{eqnarray}\lbl{Mike}
\Psi_n
 \leq  E\Big\{ P^1\Big(|\sum_{k=1}^nx_{k1}x_{k2}| > a_n\Big)^2I_{G_n}\Big\} + 2e^{-nC_{\delta}}
\end{eqnarray}
for all $n\geq 1.$ Evidently, $|x_{k1}x_{k2}|\leq C^2$, $E^1(x_{k1}x_{k2})=0$ and $E^1(x_{k1}x_{k2})^2=x_{k1}^2,$ where $E^1$ stands for the conditional expectation given $\{x_{k1},\, 1\leq k \leq n\}.$ By the Bernstein inequality (see, e.g., p.111 from Chow and Teicher (1997)),
\begin{eqnarray}\lbl{Yinmin}
P^1\Big(|\sum_{k=1}^nx_{k1}x_{k2}| > a_n\Big)^2I_{G_n}
& \leq & 4\cdot \exp\Big\{-\frac{a_n^2}{(\sum_{k=1}^nx_{k1}^2+ C^2a_n)}\Big\}I_{G_n} \nonumber\\
& \leq &  4\cdot\exp\Big\{-\frac{a_n^2}{((1+\delta)n+ C^2a_n)}\Big\}\nonumber\\
& \leq & \frac{1}{p^{t^2/(1+2\delta)}}
\end{eqnarray}
as $n$ is sufficiently large, since $a_n^2/(n(1+\delta)+ C^2a_n) \sim t^2(\log p)/(1+\delta)$ as $n\to\infty.$ Recalling (\ref{Mike}), the conclusion then follows by taking $\delta$ small enough.\ \ \ \ \ \ \ \ $\blacksquare$\\

\noindent\textbf{Proof of Lemma \ref{sweet1}}. Let $P^2$ stand for the conditional probability given $\{x_{k2},\, 1\leq k \leq n\}.$ Since $\{x_{ij};\, i\geq 1,\, j\geq 1\}$ are i.i.d., to prove the lemma, it is enough to prove
\begin{eqnarray}\lbl{uncharlie}
\Psi_n:= E\Big\{ P^2\Big(|\sum_{k=1}^nx_{k1}x_{k2}| > t_n\sqrt{n\log p}\,\Big)^2\Big\}=O\left(\frac{1}{p^{t^2-\epsilon}}\right)
\end{eqnarray}
as $n\to\infty.$ We do this only for convenience of notation.

 {\it Step 1.} For any $x>0,$ by the Markov inequality
\begin{eqnarray}\lbl{Markov}
P(\max_{1\leq k \leq n}|x_{k2}|\geq x)\leq nP(|x_{12}|\geq x) \leq Cne^{-t_0x^{\alpha}}
\end{eqnarray}
where $C=Ee^{t_0|x_{11}|^{\alpha}}<\infty.$  Second, the given condition implies that $Ee^{t|x_{11}|^{4\beta/(1+\beta)}}<\infty$ for any $t>0.$ For any $\epsilon>0,$ by (ii) of Lemma \ref{Xia}, there exists a constant $C=C_{\epsilon}>0$ such that
\begin{eqnarray}\lbl{MDP}
P\Big(\frac{|\sum_{k=1}^nx_{k2}^2-n|}{n^{(\beta + 1)/2}}\geq \epsilon \Big) \leq e^{-C_{\epsilon}n^{\beta}}
\end{eqnarray}
for each $n\geq 1.$

Set $h_n=n^{(1-\beta)/4}$, $\mu_n=E x_{ij}I(|x_{ij}| \leq h_n)$,
\begin{eqnarray}
& &  y_{ij}=x_{ij}I(|x_{ij}| \leq h_n) - E x_{ij}I(|x_{ij}| \leq h_n)\ \ \ \ \ \ \ \ \ \ \ \ \ \ \nonumber\\
& & z_{ij}=x_{ij}I(|x_{ij}| > h_n) - E x_{ij}I(|x_{ij}| > h_n)\lbl{piano}
\end{eqnarray}
for all $i\geq 1$ and $j\geq 1.$ Then, $x_{ij}=y_{ij} + z_{ij}$ for all $i, j\geq 1.$ Use the inequality $P(U+V \geq u+v)\leq P(U\geq u) + P(V\geq v)$ to obtain
\begin{eqnarray}\lbl{book}
& & P^2\Big(|\sum_{k=1}^nx_{k1}x_{k2}| > t_n\sqrt{n\log p}\,\Big)^2\nonumber\\
& \leq & 2P^2\Big(|\sum_{k=1}^ny_{k1}x_{k2}| > (t_n-\delta)\sqrt{n\log p}\,\Big)^2 +2P^2\Big(|\sum_{k=1}^nz_{k1}x_{k2}| > \delta\sqrt{n\log p}\,\Big)^2\nonumber\\
& := & 2A_n + 2B_n
\end{eqnarray}
for any $\delta>0$ small enough. Hence,
\begin{eqnarray}\lbl{Honda}
\Psi_n \leq 2EA_n + 2E B_n
\end{eqnarray}
for all $n\geq 2.$

\noindent{\it Step 2: the bound of $A_n$.} Now, if $\max_{1\leq k \leq n}|x_{k2}|\leq h_n,$ then $|y_{k1}x_{k2}|\leq 2h_n^2$ for all $k\geq 1.$ It then follows from the Bernstein inequality (see, e.g., p. 111 from Chow and Teicher (1997)) that
\begin{eqnarray*}
A_n &= & P^2\Big(|\sum_{k=1}^ny_{k1}x_{k2}| > (t_n-\delta)\sqrt{n\log p}\,\Big)^2\\
&\leq & 4\cdot\exp\Big\{-\frac{(t_n-\delta)^2n\log p}{E(y_{11}^2)\sum_{k=1}^nx_{k2}^2 + 2h_n^2(t_n-\delta)\sqrt{n\log p}}\Big\}\\
& \leq & 4\cdot\exp\Big\{-\frac{(t_n-\delta)^2n\log p}{E(y_{11}^2)(n+\epsilon n^{(\beta +1)/2}) + 2h_n^2(t_n-\delta)\sqrt{n\log p}}\Big\}
\end{eqnarray*}
 for $0<\delta < t_n$ and $\frac{|\sum_{k=1}^nx_{k2}^2-n|}{n^{(\beta + 1)/2}}< \epsilon.$   Notice $E(y_{11}^2)\to 1$ and $2h_n^2(t_n-\delta)\sqrt{n\log p}/3=o(n)$ as $n\to\infty.$ Thus,
\begin{eqnarray*}
\frac{(t_n-\delta)^2n\log p}{E(y_{11}^2)(n+\epsilon n^{(\beta +1)/2}) + 2h_n^2(t_n-\delta)\sqrt{n\log p}}\sim (t-\delta)^2\log p
\end{eqnarray*}
as $n\to\infty.$ In summary, if $\max_{1\leq k \leq n}|x_{k2}|\leq h_n$ and $\frac{|\sum_{k=1}^nx_{k2}^2-n|}{n^{(\beta + 1)/2}}\leq \epsilon,$ then for any $\delta\in (0, t/2),$
\begin{eqnarray}
A_n\leq \frac{1}{p^{t^2-2t\delta}}
\end{eqnarray}
as $n$ is sufficiently large. Therefore, for any $\epsilon>0$ small enough, take $\delta$ sufficiently small to obtain
\begin{eqnarray}\lbl{sleep}
E A_n&= &  E\Big\{ P^2\Big(|\sum_{k=1}^ny_{k1}x_{k2}| > (t_n-\delta)\sqrt{n\log p}\,\Big)^2\Big\}\nonumber\\
 & \leq & \frac{1}{p^{t^2-\epsilon}} +P(\max_{1\leq k \leq n}|x_{k2}|\geq h_n) +P\Big(\frac{|\sum_{k=1}^nx_{k2}^2-n|}{n^{(\beta + 1)/2}}\geq \epsilon \Big)\nonumber\\
 & \leq &  \frac{1}{p^{t^2-\epsilon}} + Cne^{-h_n^{\alpha}} + e^{-C_{\epsilon}n^{\beta}}=O\left(\frac{1}{p^{t^2-\epsilon}}\right)
\end{eqnarray}
as $n\to \infty,$ where the second inequality follows from (\ref{Markov}) and (\ref{MDP}), and the last identity follows from the fact that $h_n^{\alpha}=n^{\beta}$ and the assumption $\log p=o(n^{\beta}).$

\noindent{\it Step 3: the bound of $B_n$.} Recalling the definition of $z_{ij}$ and $\mu_n$ in (\ref{piano}), we have
\begin{eqnarray}
\sqrt{B_n} &= & P^2\Big(|\sum_{k=1}^nz_{k1}x_{k2}| > \delta\sqrt{n\log p}\,\Big)\nonumber\\
& \leq & P^2\Big(|\sum_{k=1}^nx_{k1}x_{k2}I\{|x_{k1}|>h_n\}| > \delta\sqrt{n\log p}/2\,\Big) + I\Big(|\sum_{k=1}^nx_{k2}| > \frac{\delta \sqrt{n\log p}}{2(e^{-n}+|\mu_n|)}\,\Big)\nonumber\\
& := & C_n + D_n.\lbl{buy}
\end{eqnarray}
Now,  by (\ref{Markov}),
\begin{eqnarray}\lbl{sell}
C_n \leq P(\max_{1\leq k\leq n}|x_{k1}|>h_n) \leq C n e^{-t_0h_n^{\alpha}}=Cn e^{-t_0n^{\beta}}.
\end{eqnarray}
Easily, $|\mu_n|\leq E |x_{11}|I(|x_{11}| > h_n)\leq e^{-t_0h_n^{\alpha}/2}E (|x_{11}|e^{t_0|x_{11}|^{\alpha}/2})=C e^{-t_0n^{\beta}/2}.$ Also, $P(|\sum_{k=1}^n\eta_k|\geq x) \leq \sum_{k=1}^nP(|\eta_k|\geq x/n)$ for any random variables $\{\eta_i\}$ and $x>0.$ We then have
\begin{eqnarray}
ED_n & = & P\Big(|\sum_{k=1}^nx_{k2}| > \frac{\delta \sqrt{n\log p}}{2(e^{-n}+|\mu_n|)}\,\Big)\nonumber\\
& \leq & n P\Big(|x_{11}| > \frac{\delta \sqrt{n\log p}}{2n(e^{-n}+|\mu_n|)} \Big)\nonumber\\
& \leq & nP\Big(|x_{11}| > e^{t_0n^{\beta}/3}\Big) \leq e^{-n}\lbl{easter}
\end{eqnarray}
as $n$ is sufficiently large, where the last inequality is from condition $Ee^{t_0|x_{11}|^{\alpha}}< \infty.$ Consequently,
\begin{eqnarray}\lbl{can}
E B_n\leq 2E(C_n^2) + 2E(D_n^2)=2E(C_n^2) + 2E(D_n)\leq e^{-Cn^{\beta}}
\end{eqnarray}
as $n$ is sufficiently large. This joint with (\ref{Honda}) and (\ref{sleep}) yields (\ref{uncharlie}). \ \ \ \ \ \ \ \ $\blacksquare$\\

\noindent\textbf{Proof of Lemma \ref{long}}. Take $\gamma = (1-\beta)/2\in [1/3, 1/2).$ Set
\begin{eqnarray}\lbl{rain}
\eta_i=\xi_iI(|\xi_i|\leq n^{\gamma}),\ \mu_n=E\eta_1\ \, \mbox{and}\ \, \sigma_n^2=Var(\eta_1),\ \ 1\leq i \leq n.
\end{eqnarray}
Since the desired result is a conclusion about $n\to\infty,$ without loss of generality, assume $\sigma_n>0$ for all $n\geq 1.$ We first claim that there exists a constant $C>0$ such that
\begin{eqnarray}\lbl{dew}
\max\Big\{|\mu_n|,\ |\sigma_n-1|,\ P(|\xi_1|>n^{\gamma}) \Big\}\leq Ce^{-n^{\beta}/C}
\end{eqnarray}
for all $n\geq 1.$ In fact, since $E\xi_1=0$ and $\alpha\gamma=\beta$,
\begin{eqnarray}\lbl{Lily}
|\mu_n| =|E\xi_1I(|\xi_1|>n^{\gamma})| \leq E|\xi_1|I(|\xi_1|>n^{\gamma})\leq E\Big(|\xi_1|e^{t_0|\xi_1|^{\alpha}/2}\Big)\cdot e^{-t_0n^{\beta}/2}
\end{eqnarray}
for all $n\geq 1.$ Note that $|\sigma_n-1|\leq |\sigma_n^2-1|=\mu_n^2 + E\xi_1^2I(|\xi_1|>n^{\gamma}),$ by the same argument as in (\ref{Lily}), we know both $|\sigma_n-1|$ and $P(|\xi_1|>n^{\gamma})$ are bounded by $Ce^{-n^{\beta}/C}$ for some $C>0.$ Then (\ref{dew}) follows.

\noindent{\it Step 1}. We prove that, for some constant $C>0,$
\begin{eqnarray}\lbl{sunshine}
\big|P\Big(\frac{S_n}{\sqrt{n\log p_n}}\geq y_n\Big)- P\Big(\frac{\sum_{i=1}^n \eta_i}{\sqrt{n\log p_n}}\geq y_n\Big)\big| \leq 2 e^{-n^{\beta}/C}
\end{eqnarray}
for all $n\geq 1.$ Observe
\begin{eqnarray}
\xi_i \equiv \eta_i\ \ \mbox{for}\ \ 1\leq i \leq n\ \ \mbox{if}\ \ \max_{1\leq i \leq n}|\xi_i| \leq n^{\gamma}.\ \ \ \ \ \ \ \ \ \ \ \ \ \ \ \ \ \ \ \ \ \ \
\end{eqnarray}
Then, by (\ref{dew}),
\begin{eqnarray}\lbl{earlier}
P\Big(\frac{S_n}{\sqrt{n\log p_n}}\geq y_n\Big)& \leq & P\Big(\frac{S_n}{\sqrt{n\log p_n}}\geq y_n, \max_{1\leq i \leq n}|\xi_i|\leq n^{\gamma}\Big) + P\Big(\bigcup_{i=1}^n\{|\xi_i| > n^{\gamma}\}\Big)\nonumber\\
& \leq & P\Big(\frac{\sum_{i=1}^n \eta_i}{\sqrt{n\log p_n}}\geq y_n\Big) + C n e^{-n^{\beta}/C}
\end{eqnarray}
for all $n\geq 1.$ Use inequality that $P(AB) \geq P(A)-P(B^c)$ for any events $A$ and $B$ to have
\begin{eqnarray*}
P\Big(\frac{S_n}{\sqrt{n\log p_n}}\geq y_n\Big)& \geq & P\Big(\frac{S_n}{\sqrt{n\log p_n}}\geq y_n, \max_{1\leq i \leq n}|\xi_i|\leq n^{\gamma}\Big)\\
& = & P\Big(\frac{\sum_{i=1}^n \eta_i}{\sqrt{n\log p_n}}\geq y_n, \max_{1\leq i \leq n}|\xi_i|\leq n^{\gamma}\Big)\\
& \geq & P\Big(\frac{\sum_{i=1}^n \eta_i}{\sqrt{n\log p_n}}\geq y_n\Big) - C n e^{-n^{\beta}/C}
\end{eqnarray*}
where in the last step the inequality $P(\max_{1\leq i \leq n}|\xi_i|> n^{\gamma})\leq C n e^{-n^{\beta}/C}$ is used as in (\ref{earlier}). This and (\ref{earlier}) concludes (\ref{sunshine}).

\noindent{\it Step 2}.  Now we prove
\begin{eqnarray}\lbl{dry}
P\Big(\frac{\sum_{i=1}^n \eta_i}{\sqrt{n\log p_n}}\geq y_n\Big) \sim \frac{e^{-x_n^2/2}}{\sqrt{2\pi} x_n}
\end{eqnarray}
as $n\to\infty,$ where
\begin{eqnarray}\lbl{land}
x_n=y_n'\sqrt{\log p_n}\ \ \mbox{and}\ \    y_n'=\frac{1}{\sigma_n}\left(y_n - \sqrt{\frac{n}{\log p_n}}\, \mu_n\right).
\end{eqnarray}
First, by (\ref{dew}),
\begin{eqnarray}\lbl{nose}
|y_n'-y_n| \leq \frac{|1-\sigma_n|}{\sigma_n}y_n + \frac{1}{\sigma_n}\cdot \sqrt{\frac{n}{\log p_n}}\,|\mu_n| \leq Ce^{-n^{\beta}/C}
\end{eqnarray}
for all $n\geq 1$ since both $\sigma_n$ and $y_n$ have limits and $p_n\to\infty.$ In particular, since $\log p_n=o(n^{\beta}),$
\begin{eqnarray}\lbl{naughty}
x_n=o(n^{\beta/2})
\end{eqnarray}
as $n\to\infty.$ Now, set
\begin{eqnarray*}
\eta_i'=\frac{\eta_i -\mu_n}{\sigma_n}
\end{eqnarray*}
for $1\leq i \leq n.$ Easily
\begin{eqnarray}\lbl{triangle}
P\Big(\frac{\sum_{i=1}^n \eta_i}{\sqrt{n\log p_n}}\geq y_n\Big) =P\Big(\frac{\sum_{i=1}^n \eta_i'}{\sqrt{n\log p_n}}\geq y_n'\Big)
\end{eqnarray}
for all $n\geq 1.$ Reviewing (\ref{rain}), for some constant $K>0,$ we have $|\eta_i'|\leq Kn^{\gamma}$ for $1\leq i \leq n.$ Take $c_n=Kn^{\gamma-1/2}.$ Recalling $x_n$ in (\ref{land}). It is easy to check that
\begin{eqnarray*}
s_n:=\Big(\sum_{i=1}^nE\eta_i'^2\Big)^{1/2}=\sqrt{n},\ \varrho_n:=\sum_{i=1}^nE|\eta_i'|^3\sim n C,\ \     |\eta_i'| \leq c_ns_n\ \mbox{and}\ 0<c_n \leq  1
\end{eqnarray*}
as $n$ is sufficiently large. Recall $\gamma = (1-\beta)/2,$ it is easy to see from (\ref{naughty}) that
\begin{eqnarray*}
0<x_n<\frac{1}{18c_n}
\end{eqnarray*}
for $n$ large enough. Now, let $\gamma(x)$ be as in Lemma \ref{shao}, since $\beta \leq 1/3,$ by the lemma and (\ref{naughty}),
\begin{eqnarray*}
\Big|\gamma\big(\frac{x_n}{s_n}\big)\Big| \leq \frac{2x_n^3\varrho_n}{s_n^3}=o\left(n^{\frac{3\beta}{2} -\frac{1}{2}}\right)\to 0\ \ \mbox{and}\ \ \frac{(1+x_n)\varrho_n}{s_n^3}=O(n^{(\beta-1)/2})\to 0
\end{eqnarray*}
as $n\to\infty.$ By (\ref{land}) and (\ref{nose}), $x_ns_n= y_n'\sqrt{n\log p_n}$ and $x_n\to\infty$ as $n\to\infty.$ Use Lemma \ref{shao} and the fact $1-\Phi(t)=\frac{1}{\sqrt{2\pi} t}e^{-t^2/2}$ as $t\to +\infty$ to obtain
\begin{eqnarray}
P\Big(\frac{\sum_{i=1}^n \eta_i'}{\sqrt{n\log p_n}}\geq y_n'\Big) =P\Big(\sum_{i=1}^n \eta_i' \geq x_ns_n\Big)\sim 1-\Phi(x_n)\sim \frac{e^{-x_n^2/2}}{\sqrt{2\pi} x_n}
\end{eqnarray}
as $n\to\infty.$ This and (\ref{triangle}) conclude (\ref{dry}).

\noindent{\it Step 3}. Now we show
\begin{eqnarray}\lbl{detail}
\frac{e^{-x_n^2/2}}{\sqrt{2\pi x_n}} \sim \frac{p_n^{-y_n^2/2}(\log p_n)^{-1/2}}{\sqrt{2\pi} y}:=\omega_n
\end{eqnarray}
as $n\to\infty.$ Since $y_n \to y$ and $\sigma_n\to 1,$ we know from (\ref{nose}) that
\begin{eqnarray}\lbl{fraction}
\sqrt{2\pi} x_n =\sqrt{2\pi} y_n'(\log p_n)^{1/2} \sim \sqrt{2\pi} y\,(\log p_n)^{1/2}
\end{eqnarray}
as $n\to\infty.$ Further, by (\ref{land}),
\begin{eqnarray}\lbl{hair}
\frac{e^{-x_n^2/2}}{p_n^{-y_n^2/2}} = \exp\Big\{-\frac{x_n^2}{2} + \frac{y_n^2}{2}\log p_n\Big\}= \exp\Big\{\frac{1}{2}\Big(y_n^2-y_n'^2\Big)\log p_n\Big\}.
\end{eqnarray}
Since $y_n\to y,$ by (\ref{nose}), both $\{y_n\}$ and $\{y_n'\}$ are bounded. It follows from (\ref{nose}) again that $|y_n^2-y_n'^2| \leq C |y_n-y_n'| =O(e^{-n^{\beta}/C})$ as $n\to\infty.$ With assumption $\log p_n=o(n^{\beta})$ we get $e^{-x_n^2/2} \sim p_n^{-y_n^2/2}$ as $n\to\infty,$ which combining with (\ref{fraction}) yields (\ref{detail}).

Finally, we compare the right hand sides of (\ref{sunshine}) and (\ref{detail}). Choose $C'>\max\{y_n^2;\, n\geq 1\},$ since $\log p_n =o(n^{\beta}),$ recall $\omega_n$ in (\ref{detail}),
\begin{eqnarray*}
\frac{2e^{-n^{\beta}/C}}{\omega_n} & = & 2\sqrt{2\pi}\, y\,(\log p_n)^{1/2}p_n^{y_n^2/2} e^{-n^{\beta}/C}\\
& = & O\left(n^{\beta/2}\cdot\exp\Big\{C'\log p_n-\frac{n^{\beta}}{C} \Big\}\right)\\
&=& O\left(n^{\beta/2}\cdot\exp\Big\{-\frac{n^{\beta}}{2C}\Big\}\right)\to 0
\end{eqnarray*}
as $n\to\infty$ for any constant $C>0.$ This fact joint with (\ref{sunshine}), (\ref{dry}) and (\ref{detail}) proves the lemma.\ \ \ \ \ \ \ \ $\blacksquare$\\

\noindent\textbf{Proof of Lemma \ref{webster}}. For any Borel set $A \subset \mathbb{R},$ set $P_2(A)=P(A|u_{k1}, u_{k3},\, 1\leq k\leq n),$ the conditional probability of $A$ with respect to $u_{k1}, u_{k3},\, 1\leq k \leq n.$ Observe from the expression of $\Sigma_4$ that three sets of random variables $\{u_{k1}, u_{k3};\, 1\leq k \leq n\}$, $\{u_{k2};\, 1\leq k \leq n\}$ and $\{u_{k4};\, 1\leq k \leq n\}$ are independent. Then
\begin{eqnarray*}
& & P\Big(|\sum_{k=1}^n u_{k1}u_{k2}|> a_n,\, |\sum_{k=1}^n u_{k3}u_{k4}|> a_n\Big)\\
& = & E\Big\{P_2\Big(|\sum_{k=1}^n u_{k1}u_{k2}|> a_n\Big)P_2\Big(|\sum_{k=1}^n u_{k3}u_{k4}|> a_n\Big)\Big\}\\
& \leq & \Big\{E\, P_2\Big(|\sum_{k=1}^n u_{k1}u_{k2}|> a_n\Big)^2\Big\}^{1/2}\cdot \Big\{E\, P_2\Big(|\sum_{k=1}^n u_{k3}u_{k4}|> a_n\Big)^2\Big\}^{1/2}
\end{eqnarray*}
by the Cauchy-Schwartz inequality. Use the same independence again
\begin{eqnarray}
& & P_2\Big(|\sum_{k=1}^n u_{k1}u_{k2}|> a_n\Big)=P\Big(|\sum_{k=1}^n u_{k1}u_{k2}|> a_n \Big|u_{k1},\, 1\leq k \leq n\Big);\lbl{cautious1}\\
& & P_2\Big(|\sum_{k=1}^n u_{k3}u_{k4}|> a_n\Big)=P\Big(|\sum_{k=1}^n u_{k3}u_{k4}|> a_n\Big|u_{k3},\, 1\leq k \leq n\Big)\lbl{cautious2}.
\end{eqnarray}
These can be also seen from  Proposition 27 in Fristedt and Gray (1997). It follows that
\begin{eqnarray*}
& & \sup_{|r|\leq 1}P\Big(|\sum_{k=1}^n u_{k1}u_{k2}|> a_n,\, |\sum_{k=1}^n u_{k3}u_{k4}|> a_n\Big)\\
& \leq & E \Big\{ P\Big(|\sum_{k=1}^n u_{k1}u_{k2}|> a_n \Big|u_{11},\cdots, u_{n1}\Big)^2\Big\}.
\end{eqnarray*}
Since $\{u_{k1};\, 1\leq k \leq n\}$ and $\{u_{k2};\, 1\leq k \leq n\}$ are independent, and $t_n:=a_n/\sqrt{n\log p}\to t=2,$ taking $\alpha=2$ in Lemma \ref{sweet1}, we obtain the desired conclusion from the lemma.\ \ \ \ \ \ \ \ $\blacksquare$\\

\noindent\textbf{Proof of Lemma \ref{yahoo}}. Since $\Sigma_4$ is always non-negative definite, the determinant of the first $3\times 3$ minor of $\Sigma_4$ is non-negative: $1-r_1^2-r_2^2 \geq 0.$  Let $r_3=\sqrt{1-r_1^2 -r_2^2}$ and $\{u_{k5};\, 1\leq k \leq n\}$ be i.i.d. standard normals which are independent of $\{u_{ki};\, 1\leq i\leq 4;\, 1\leq k \leq n\}.$ Then,
\begin{eqnarray*}
(u_{11}, u_{12}, u_{13}, u_{14})\overset{d}{=}(u_{11}, u_{12}, r_1u_{11} + r_2u_{12} + r_3u_{15}, u_{14}).
\end{eqnarray*}
Define $Z_{ij}=|\sum_{k=1}^n u_{ki}u_{kj}|$ for $1\leq i, j\leq 5$ and $r_5=r_3.$ By the Cauchy-Schwartz inequality,
\begin{eqnarray*}
 |\sum_{k=1}^n (r_1u_{k1} + r_2u_{k2} + r_3u_{k5}) u_{k4}| & \leq & \sum_{i\in \{1,2,5\}}|r_i|\cdot |\sum_{k=1}^n u_{ki}u_{k4}|\\
 &\leq & \Big(r_1^2 + r_2^2 + r_3^2\Big)^{1/2}
 \Big(Z_{14}^2 + Z_{24}^2 + Z_{54}^2\Big)^{1/2}\\
 & \leq & \sqrt{3}\cdot \max\{Z_{14}, Z_{24}, Z_{54}\}.
\end{eqnarray*}
It follows from the above two facts that
\begin{eqnarray}
& & P\Big(|\sum_{k=1}^n u_{k1}u_{k2}|> a_n,\, |\sum_{k=1}^n u_{k3}u_{k4}|> a_n\Big)\nonumber\\
& \leq  & P\Big(Z_{12}> a_n,\, \max\{Z_{14}, Z_{24}, Z_{54}\}> \frac{a_n}{\sqrt{3}}\Big)\nonumber\\
& \leq & \sum_{i\in \{1,2,5\}}P\Big(Z_{12}> a_n,\, Z_{i4}> \frac{a_n}{\sqrt{3}}\Big) \nonumber\\
& = & 2P\Big(Z_{12}> a_n,\, Z_{14}> \frac{a_n}{\sqrt{3}}\Big) + P\Big(Z_{12}> a_n\Big)\cdot P\Big(Z_{54}> \frac{a_n}{\sqrt{3}}\Big) \lbl{least}
\end{eqnarray}
by symmetry and independence.  For any Borel set $A \subset \mathbb{R},$ set $P^1(A)=P(A|u_{k1},\, 1\leq k\leq n),$ the conditional probability of $A$ with respect to $u_{k1}, \, 1\leq k \leq n.$ For any $s>0,$ from the fact that $\{u_{k1}\}, \{u_{k2}\}$ and $\{u_{k4}\}$ are independent, we see that
\begin{eqnarray*}
P\Big(Z_{12}> a_n,\, Z_{14}> s a_n\Big)
& = & E\Big(P^1(Z_{12}> a_n)\cdot P^1(Z_{14}> sa_n)\Big)\\
& \leq &  \Big\{E\, P^1(Z_{12}> a_n)^2\Big\}^{1/2}\cdot \Big\{E\, P^1(Z_{14}> s a_n)^2\Big\}^{1/2}
\end{eqnarray*}
by the Cauchy-Schwartz inequality. Taking $t_n:=a_n/\sqrt{n\log p}\to t=2$ and $t_n:=s a_n/\sqrt{n\log p}\to t=2s$ in Lemma \ref{sweet1}, respectively, we get
\begin{eqnarray*}
E\, P^1(Z_{12}> a_n)^2 =O\Big(p^{-4+\epsilon}\Big)\ \ \mbox{and}\ \ E P^1(Z_{14}> s a_n)^2 =O\Big(p^{-4s^2+\epsilon}\Big)
\end{eqnarray*}
as $n \to \infty$ for any $\epsilon>0.$ This implies that, for any $s>0$ and $\epsilon>0,$
\begin{eqnarray}\lbl{teamates}
P\Big(Z_{12}> a_n,\, Z_{14}> s a_n\Big) \leq O\Big(p^{-2-2s^2+\epsilon}\Big)
\end{eqnarray}
as $n \to \infty.$ In particular,
\begin{eqnarray}\lbl{bai}
P\Big(Z_{12}> a_n,\, Z_{14}> \frac{a_n}{\sqrt{3}}\Big) \leq O\Big(p^{-\frac{8}{3}+\epsilon}\Big)
\end{eqnarray}
as $n \to \infty$ for any $\epsilon>0.$

Now we bound the last term in (\ref{least}). Note that $|u_{11}u_{12}|\leq (u_{11}^2 + u_{12}^2)/2,$ it follows that $Ee^{|u_{11}u_{12}|/2}<\infty$ by independence and $E\exp(N(0,1)^2/4)<\infty.$ Since $\{u_{k1}, u_{k2};\, 1\leq k \leq n\}$ are i.i.d. with mean zero and variance one, and $y_n:=a_n/\sqrt{n\log p}\to 2$ as $n\to\infty,$ taking $\alpha=1$ in Lemma \ref{long}, we get
\begin{eqnarray}
 P\Big(Z_{12}> a_n\Big)
 &= & P\Big(\frac{1}{\sqrt{n\log p}}|\sum_{k=1}^nu_{k1}u_{k2}|> \frac{a_n}{\sqrt{n\log p}}\Big)\nonumber\\
& \sim & 2\cdot \frac{p^{-y_n^2/2}(\log p)^{-1/2}}{2\sqrt{2\pi}}\sim \frac{e^{-y/2}}{\sqrt{2\pi}}\cdot\frac{1}{p^2}\lbl{Gooden10}
\end{eqnarray}
as $n\to\infty.$ Similarly, for any $t>0$,
\begin{eqnarray}\lbl{beat}
 P\Big(Z_{12}> ta_n\Big)=O\Big(p^{-{2t^2+\epsilon}}\Big)
\end{eqnarray}
as $n\to\infty$ (this can also be derived from (i) of Lemma \ref{Xia}). In particular,
\begin{eqnarray}\lbl{omit}
 P\Big(Z_{54}> \frac{a_n}{\sqrt{3}}\Big) & = & P\Big(Z_{12}> \frac{a_n}{\sqrt{3}}\Big)=O\Big(p^{-\frac{2}{3}+\epsilon}\Big)
\end{eqnarray}
as $n\to\infty$ for any $\epsilon>0.$ Combining (\ref{Gooden10}) and (\ref{omit}), we know that the last term in (\ref{least}) is bounded by $O(p^{-\frac{8}{3}+\epsilon})$ as $n\to\infty$ for any $\epsilon>0.$ This together with (\ref{least}) and (\ref{bai}) concludes the lemma. \ \ \ \ \ \ \ \ $\blacksquare$\\

\noindent\textbf{Proof of Lemma \ref{yahoo1}}. Fix $\delta\in (0, 1).$ Take independent standard normals $\{u_{k5}, u_{k6};\, 1\leq k \leq n\}$ that are also independent of $\{u_{ki};\, 1\leq i\leq 4;\, 1\leq k \leq n\}.$ Then, since $\{u_{k1}, u_{k2}, u_{k5}, u_{k6};\, 1\leq k \leq n\}$ are i.i.d. standard normals, by checking covariance matrix $\Sigma_4$, we know
\begin{eqnarray}\lbl{toy}
(u_{11}, u_{12}, u_{13}, u_{14})\overset{d}{=}(u_{11}, u_{12}, r_1u_{11} + r_1'u_{15}, r_2u_{12} + r_2'u_{16})
\end{eqnarray}
where $r_1'=\sqrt{1-r_1^2}$ and $r_2'=\sqrt{1-r_2^2}.$ Define $Z_{ij}=|\sum_{k=1}^n u_{ki}u_{kj}|$ for $1\leq i, j\leq 6.$ Then
\begin{eqnarray}
& & |\sum_{k=1}^n (r_1u_{k1} + r_1'u_{k5})(r_2u_{k2} + r_2'u_{k6})|\nonumber\\
& \leq &  |r_1r_2| Z_{12} + |r_1r_2'| Z_{16} + |r_1'r_2| Z_{25} +|r_1'r_2'| Z_{56}\nonumber\\
& \leq & (1-\delta)^2 Z_{12} + 3\max\{Z_{16}, Z_{25}, Z_{56}\}\lbl{Ruthy}
\end{eqnarray}
for all $|r_1|,\, |r_2|\leq 1-\delta.$ Let $\alpha=(1+ (1-\delta)^2)/2$, $\beta=\alpha/(1-\delta)^2$ and $\gamma=(1-\alpha)/3.$ Then
\begin{eqnarray}\lbl{shot}
\beta>1\ \  \mbox{and}\ \ \gamma>0.
\end{eqnarray}
Easily, if $Z_{12}\leq \beta a_n$, $\max\{Z_{16}, Z_{25}, Z_{56}\}\leq \gamma a_n,$ then from (\ref{Ruthy}) we know that the left hand side of (\ref{Ruthy}) is controlled by $a_n.$ Consequently, by (\ref{toy}) and the i.i.d. property,
\begin{eqnarray}
P(Z_{12}>a_n,\, Z_{34}>a_n) &= & P\Big(Z_{12}>a_n,\, |\sum_{k=1}^n (r_1u_{k1} + r_1'u_{k5})(r_2u_{k2} + r_2'u_{k6})|>a_n\Big)\nonumber\\
& \leq & P(Z_{12}>a_n, Z_{12}>\beta a_n) +  \sum_{i\in \{1,2,5\}}P(Z_{12}>a_n, Z_{i6}>\gamma a_n)\nonumber\\
& = & P(Z_{12}>\beta a_n) + 2P(Z_{12}>a_n, Z_{16}>\gamma a_n)\nonumber\\
& & \ \ \ \ \ \ \ \ \ \ \ \ \ \ \ \ \ \ \   +\, P(Z_{12}>a_n)\cdot P(Z_{56}>\gamma a_n)\lbl{cough}
\end{eqnarray}
where ``$2P(Z_{12}>a_n, Z_{16}>\gamma a_n)$" comes from the fact $(Z_{12}, Z_{16})\overset{d}{=}(Z_{12}, Z_{26}).$ Keep in mind that $(Z_{12}, Z_{16})\overset{d}{=}(Z_{12}, Z_{14})$ and $Z_{56}\overset{d}{=}Z_{12}.$ Recall (\ref{shot}), applying (\ref{teamates}) and (\ref{beat}) to the three terms in the sum on the right hand side of (\ref{cough}), we conclude (\ref{demand}).\ \ \ \ \ \ \ $\blacksquare$\\

\noindent\textbf{Proof of Lemma \ref{fertilizer}}. Reviewing notation $\Omega_3$ defined below (\ref{fire}), the current case is that $d_1\leq d_3\leq d_2 \leq d_4$ with $d=(d_1, d_2)$ and $d'=(d_3, d_4).$ Of course, by definition, $d_1< d_2$ and $d_3< d_4.$ To save notation, define the ``neighborhood" of $d_i$ as follows:
\begin{eqnarray}
N_i=\Big\{d\in \{1,\cdots, p\};\, |d-d_i|<\tau\Big\}\ \ \ \ \ \ \ \ \ \ \ \ \ \ \ \ \ \ \ \
\end{eqnarray}
for $i=1,2,3,4.$

Given $d_1< d_2,$ there are two possibilities for $d_4$: (a) $d_4-d_2> \tau$ and (b) $0 \leq d_4-d_2 \leq \tau.$ There are four possibilities for $d_3$: (A) $ d_3\in N_2\backslash N_1$; (B) $d_3\in N_1\backslash N_2$; (C) $d_3\in N_1\cap N_2$; (D) $ d_3\notin N_1\cup N_2$. There are eight combinations for the locations of $(d_3, d_4)$ in total. However, by (\ref{helicopter}) the combination (a) \& (D) is excluded. Our analysis next will exhaust all of the seven possibilities.

\noindent\textbf{Case (a) \& (A)}. Let $\Omega_{a, A}$ be the subset of $(d, d')\in \Omega_3$ satisfying restrictions (a) and (A), and others such as $\Omega_{b, C}$ are  similarly defined. Thus,
\begin{eqnarray}\lbl{nightfall}
 \sum_{(d, d')\in \Omega_3} P(Z_{d} >a_n, Z_{d'} > a_n)\leq  \sum_{\theta, \Theta}\sum _{(d, d')\in \Omega_{\theta, \Theta}} P(Z_{d} >a_n, Z_{d'} > a_n)
\end{eqnarray}
where $\theta$ runs over set $\{a, b\}$ and $\Theta$ runs over set $\{A, B,C, D\}$ but $(\theta, \Theta) \ne (a, D).$

Easily, $|\Omega_{a, A}|\leq \tau p^3$ and the covariance matrix of $(w_{d_2}, w_{d_1}, w_{d_3}, w_{d_4})$ (see (\ref{coordinate})) is
\begin{eqnarray*}
\begin{pmatrix}
1& 0 & \gamma & 0\\
0 & 1 & 0 & 0\\
\gamma & 0 & 1 & 0\\
0 & 0 & 0 & 1
\end{pmatrix}
,\ \ |\gamma| \leq 1.
\end{eqnarray*}
Take $\epsilon=1/2$ in Lemma \ref{webster} to have $P(Z_{d} > a_n, Z_{d'} > a_n)\equiv \rho_n=o(p^{-7/2})$ for all $(d, d')\in\Omega_{a, A}.$ Thus
\begin{eqnarray}\lbl{aA}
\sum_{(d, d')\in R} P(Z_{d} > a_n, Z_{d'} > a_n)=|R|\cdot\rho_n\to 0
\end{eqnarray}
as $n\to\infty$ for $R=\Omega_{a, A}.$

\noindent\textbf{Case (a) \& (B)}. Notice $|\Omega_{a, B}|\leq \tau p^3$ and the covariance matrix of $(w_{d_1}, w_{d_2}, w_{d_3}, w_{d_4})$ is the same as that in Lemma \ref{webster}. By the lemma we then have  (\ref{aA})  for  $R=\Omega_{a, B}.$

\noindent\textbf{Case (a) \& (C)}. Notice $|\Omega_{a, C}|\leq \tau^2p^2$ and the covariance matrix of $(w_{d_1}, w_{d_2}, w_{d_3}, w_{d_4})$ is the same as that in Lemma \ref{yahoo}. By the lemma, we know  (\ref{aA}) holds for  $R=\Omega_{a, C}.$

\noindent\textbf{Case (b) \& (A)}. In this case, $|\Omega_{b, A}|\leq \tau^2p^2$ and the covariance matrix of $(w_{d_3}, w_{d_4}, w_{d_2}, w_{d_1})$ is the same as that in Lemma \ref{yahoo}. By the lemma and using the fact that
\begin{eqnarray*}
P(Z_{d} > a_n, Z_{d'} > a_n)=P(Z_{(d_3, d_4)} > a_n, Z_{(d_2, d_1)} > a_n)
\end{eqnarray*}
we see  (\ref{aA}) holds with  $R=\Omega_{b, A}.$

\noindent\textbf{Case (b) \& (B)}. In this case, $|\Omega_{b, B}|\leq \tau^2p^2$ and the covariance matrix of $(w_{d_1}, w_{d_2}, w_{d_3}, w_{d_4})$ is the same as that in Lemma \ref{yahoo1}. By the lemma, we know  (\ref{aA}) holds for  $R=\Omega_{b, B}.$


\noindent\textbf{Case (b) \& (C)}. We assign positions for $d_1, d_3, d_2, d_4$ step by step: there are at most $p$ positions for $d_1$ and at most $k$ positions for each of $d_3, d_2$ and $d_4.$ Thus, $|\Omega_{b, C}|\leq \tau^3p.$ By (\ref{Gooden10}),
\begin{eqnarray*}
P(Z_{d} > a_n, Z_{d'} > a_n) \leq P(Z_{d} > a_n)= P\Big(|\sum_{i=1}^n\xi_i\eta_i|>a_n\Big)=O\Big(\frac{1}{p^2}\Big)
\end{eqnarray*}
as $n\to\infty,$ where $\{\xi_i, \eta_i;\, i\geq 1\}$ are i.i.d. standard normals. Therefore, (\ref{aA}) holds with  $R=\Omega_{b, C}.$

\noindent\textbf{Case (b) \& (D)}. In this case, $|\Omega_{b, C}|\leq \tau p^3$ and the covariance matrix of $(w_{d_4}, w_{d_3}, w_{d_2}, w_{d_1})$ is the same as that in Lemma \ref{webster}. By the lemma and noting the fact that
\begin{eqnarray*}
P(Z_{d} > a_n, Z_{d'} > a_n)=P(Z_{(d_4, d_3)} > a_n, Z_{(d_2, d_1)} > a_n)
\end{eqnarray*}
we see  (\ref{aA}) holds with  $R=\Omega_{b, D}.$

We obtain (\ref{Suzuki}) by combining (\ref{aA}) for all the cases considered above with (\ref{nightfall}). \ \ \ \ \ \ \ $\blacksquare$


\end{document}